\documentclass[11pt]{article}
\usepackage{makeidx}         % allows index generation
\usepackage{graphicx}        % standard LaTeX graphics tool
\usepackage{psfrag}                             % when including figure files
\usepackage{multicol}        % used for the two-column index
\usepackage{amsmath}
\usepackage{amssymb}
\usepackage{amsfonts}
\usepackage{latexsym}
\usepackage[T1]{fontenc}
\usepackage{subfigure}                      % figures and subfigures
\usepackage{float}
\usepackage{mathrsfs}
\newtheorem{theorem}{Theorem}
\newtheorem{proposition}{Proposition}
\newtheorem{definition}{Definition}
\newtheorem{remark}{Remark}
\newtheorem{lemma}{Lemma}
\newtheorem{corollary}{Corollary}
\newtheorem{example}{Example}

\begin{document}

\title{A survey on the boundary behavior of the double layer potential in Schauder spaces in the frame of an abstract approach}

\index{boundary behavior of the double layer potential}

%\titlerunning{Boundary behavior of the double layer potential in H\"{o}lder spaces}

\author{  
Massimo Lanza de Cristoforis
\\
Dipartimento di Matematica `Tullio Levi-Civita', 
\\
Universit\`a degli Studi di Padova, 
\\
Via Trieste 63, Padova 35121, 
Italy. 
\\
E-mail: mldc@math.unipd.it   }
 \date{\ }
%\authorrunning{M.~Lanza de Cristoforis}

%\institute{M.~Lanza de Cristoforis\at Dipartimento di Matematica 
%`Tullio Levi-Civita', 
%Universit\`{a} di Padova,\hfill\break
%\email{mldc@math.unipd.it} }

\maketitle

%\index{!}
%\index{}

\bgroup

% The following are the author's own macros. Use your own ONLY if you are
% conversant with LaTeX and know what you are doing. There must be no clash
% between your macros and the publisher's. To avoid this, use the construction
% with \bgroup (see immediately above) and \par\egroup (see at the very end
% of this specimen paper). This will make your macros local.

\allowdisplaybreaks

\abstract{We provide a summary of  the continuity properties  of the boundary integral operator corresponding to the double layer potential associated to the fundamental solution of a {\em nonhomogeneous} second order elliptic differential operator with constant coefficients
   in H\"{o}lder and Schauder spaces  on the boundary of a bounded open   subset of ${\mathbb{R}}^n$.
   The purpose is two-fold. On one hand we try present in a single paper all the known continuity results on the topic   with the best known exponents
    in a H\"{o}lder and Schauder space setting  and on the other hand we show that many of the properties we present can be deduced by applying results that hold in an abstract setting of metric spaces with a measure that satisfies certain growth conditions that include   non-doubling measures as  in a series of papers by   Garc\'{\i}a-Cuerva and Gatto   in the frame of H\"{o}lder spaces and later by the author.}

\section{Introduction}  In this paper, we consider the double layer potential associated to the fundamental solution of a second order differential operator with constant coefficients. Throughout the paper, we assume that
\[
n\in {\mathbb{N}}\setminus\{0,1\}\,,
\]
where ${\mathbb{N}}$ denotes the set of natural numbers including $0$. Let $\alpha\in]0,1]$, $m\in {\mathbb{N}} $. Let $\Omega$ be a bounded open subset of ${\mathbb{R}}^{n}$ of class $C^{m,\alpha}$. Here we understand that $C^{m,0}\equiv C^m$.    For the notation and standard properties of the (generalized) Schauder spaces and sets of class $C^{m,\alpha}$ we refer for example  to Dondi and the author  \cite[\S 2]{DoLa17} and to the reference \cite[\S 2.11, 2.13]{DaLaMu21} of Dalla Riva, the author and Musolino. 

Let $\nu \equiv (\nu_{l})_{l=1,\dots,n}$ denote the external unit normal to $\partial\Omega$. If $\Omega$ is an open Lipschitz set $\nu$ is known to exist only for almost all points of $\partial\Omega$. Let $N_{2}$ denote the number of multi-indexes $\gamma\in {\mathbb{N}}^{n}$ with $|\gamma|\leq 2$. For each 
\begin{equation}
\label{introd0}
{\mathbf{a}}\equiv (a_{\gamma})_{|\gamma|\leq 2}\in {\mathbb{C}}^{N_{2}}\,, 
\end{equation}
we set 
\[
a^{(2)}\equiv (a_{lj} )_{l,j=1,\dots,n}\qquad
a^{(1)}\equiv (a_{j})_{j=1,\dots,n}\qquad
a\equiv a_{0}\,.
\]
with $a_{lj} \equiv 2^{-1}a_{e_{l}+e_{j}}$ for $j\neq l$, $a_{jj} \equiv
 a_{e_{j}+e_{j}}$,
and $a_{j}\equiv a_{e_{j}}$, where $\{e_{j}:\,j=1,\dots,n\}$  is the canonical basis of ${\mathbb{R}}^{n}$. We note that the matrix $a^{(2)}$ is symmetric. 
Then we assume that 
  ${\mathbf{a}}\in  {\mathbb{C}}^{N_{2}}$ satisfies the following ellipticity assumption
\begin{equation}
\label{ellip}
\inf_{
\xi\in {\mathbb{R}}^{n}, |\xi|=1
}{\mathrm{Re}}\,\left\{
 \sum_{|\gamma|=2}a_{\gamma}\xi^{\gamma}\right\} >0\,,
\end{equation}
and we consider  the case in which
\begin{equation}
\label{symr}
a_{lj} \in {\mathbb{R}}\qquad\forall  l,j=1,\dots,n\,.
\end{equation}
Then we introduce the operators
\begin{eqnarray*}
%\label{introd1}
P[{\mathbf{a}},D]u&\equiv&\sum_{l,j=1}^{n}\partial_{x_{l}}(a_{lj}\partial_{x_{j}}u)
+
\sum_{l=1}^{n}a_{l}\partial_{x_{l}}u+au\,,
\\
%\label{introd2}
B_{\Omega}^{*}v&\equiv&\sum_{l,j=1}^{n} \overline{a}_{jl}\nu_{l}\partial_{x_{j}}v
-\sum_{l=1}^{n}\nu_{l}\overline{a}_{l}v\,,
\end{eqnarray*}
for all $u,v\in C^{2}(\overline{\Omega})$, and a fundamental solution $S_{{\mathbf{a}} }$ of $P[{\mathbf{a}},D]$, and the  boundary integral operator corresponding to the  double layer potential 
\begin{eqnarray}
\label{introd3}
\lefteqn{
 W_\Omega[{\mathbf{a}},S_{{\mathbf{a}}}   ,\mu](x) \equiv 
{\mathrm{p.v.}}\int_{\partial\Omega}\mu (y)\overline{B^{*}_{\Omega,y}}\left(S_{{\mathbf{a}}}(x-y)\right)
\,d\sigma_{y}
}
\\  \nonumber
&&
\qquad
=-{\mathrm{p.v.}}\int_{\partial\Omega}\mu(y)\sum_{l,j=1}^{n} a_{jl}\nu_{l}(y)\frac{\partial S_{ {\mathbf{a}} } }{\partial x_{j}}(x-y)\,d\sigma_{y}
\\  \nonumber
&&
\qquad\quad
-\int_{\partial\Omega}\mu(y)\sum_{l=1}^{n}\nu_{l}(y)a_{l}
S_{ {\mathbf{a}} }(x-y)\,d\sigma_{y}  
\end{eqnarray}
for almost all $x\in \partial\Omega$,
where the density or moment $\mu$ is a function  from $\partial\Omega$ to ${\mathbb{C}}$.  Here the subscript $y$ of $\overline{B^{*}_{\Omega,y}}$ means that we are taking $y$ as variable of the differential operator $\overline{B^{*}_{\Omega,y}}$, $d\sigma$ is the ordinary $(n-1)$-dimensional measure, and `${\mathrm{p.v.}}$' denotes the principal value of the integral. Thus the kernel of the boundary integral operator corresponding to the  double layer potential 
  is the following
\begin{eqnarray}\label{eq:tgdlgen}
\lefteqn{
K_{\Omega,{\mathbf{a}}}(x,y)\equiv \overline{B^{*}_{\Omega,y}}\left(S_{{\mathbf{a}}}(x-y)\right) 
}
\\ \nonumber
&&\qquad\qquad
\equiv - \sum_{l,j=1}^{n} a_{jl}\nu_{l}(y)\frac{\partial S_{ {\mathbf{a}} } }{\partial x_{j}}(x-y) 
 - \sum_{l=1}^{n}\nu_{l}(y)a_{l}
S_{ {\mathbf{a}} }(x-y) 
\end{eqnarray}
for all $x\in\partial\Omega$ and for almost all $y\in\partial\Omega$ with $x\neq y$
(cf.~(\ref{introd3})). The role of the double layer potential in the solution of boundary value problems for the operator $P[{\mathbf{a}},D]$ is well known (cf. \textit{e.g.}, 
G\"{u}nter~\cite{Gu67}, Kupradze,  Gegelia,  Basheleishvili and 
 Burchuladze~\cite{KuGeBaBu79}, Mikhlin \cite{Mik70},  Mikhlin and  Pr\"{o}ssdorf \cite{MikPr86}, Buchukuri,  Chkadua,  Duduchava, and  Natroshvili \cite{BuChDuNa12}.)

Here we provide a summary of the continuity  properties of the boundary operator $W_\Omega[{\mathbf{a}},S_{{\mathbf{a}}},\cdot]$ (the so-called Neumann-Poincar\'{e} operator in case $P[{\mathbf{a}},D]$ is the Laplace operator) in the frame of H\"{o}lder and Schauder spaces. 

Also, we give references to proofs where a consistent part of the arguments are based on results for integral operators that hold in a  metric space with a measure that satisfies certain growth conditions that include   non-doubling measures as  in a series of papers by   Garc\'{\i}a-Cuerva and Gatto \cite{GaGa04}, \cite{GaGa05}, Gatto \cite{Gat09} in the frame of H\"{o}lder spaces and that have been further developed  in \cite{La22a}.  We now briefly present such abstract setting, that this paper shows to have several applications
(see also \cite{La22b}, \cite{La22d}, \cite{La23b}, \cite{La23c}). Let 
$
(M,d)$ be a metric space and let 
$X$, $Y$ be subsets of $M$. 
\begin{eqnarray} \nonumber
&&\text{Let}\ {\mathcal{N}} \ \text{ be a $\sigma$-algebra of parts of}\  Y  \,, {\mathcal{B}}_Y\subseteq  {\mathcal{N}}\,.
\\   \label{eq:nu}
&&\text{Let}\ \nu\  \text{be   measure on}\  {\mathcal{N}} \,.
\\  \nonumber
&&\text{Let}\ \nu(B(x,r)\cap Y)<+\infty\qquad\forall (x,r)\in X\times ]0,+\infty[\,,
\end{eqnarray}
where  ${\mathcal{B}}_Y$ denotes the $\sigma$-algebra of the Borel subsets of $Y$ and 
\[
B(\xi,r)\equiv \left\{\eta\in M:\, d(\xi,\eta)<r\right\}\qquad\forall (\xi,r)\in M\times ]0,+\infty[\,.
\]
 We assume that  $\upsilon_Y\in ]0,+\infty[$ and we consider two types of assumptions on $\nu$. 
The first assumption is  that
$Y$ is   upper $\upsilon_Y$-Ahlfors regular with respect to $X$, \textit{i.e.}, that
\begin{eqnarray} \nonumber
&&\text{there\ exist}\ r_{X,Y,\upsilon_Y}\in]0,+\infty]\,,\ c_{X,Y,\upsilon_Y}\in]0,+\infty[\ \text{such\ that}
\\ \nonumber
&&\nu( B(x,r)\cap Y )\leq c_{X,Y,\upsilon_Y} r^{\upsilon_Y}  
\\  \label{defn:uareg1}
&&\text{for\ all}\ x\in X\ \text{and}\  r \in]0,r_{X,Y,\upsilon_Y}[
 \,. 
\end{eqnarray}
In case $X=Y$, we just say that $Y$ is   upper $\upsilon_Y$-Ahlfors regular  and this is the assumption that has been considered by 
Garc\'{\i}a-Cuerva and Gatto \cite{GaGa04}, \cite{GaGa05}, Gatto \cite{Gat06}, \cite{Gat09}   in case $X=Y=M$. See also  Edmunds,  Kokilashvili and   Meskhi~\cite[Chap.~6]{EdKoMe02} in the frame of Lebsgue spaces. 

An interesting feature of condition (\ref{defn:uareg1}) is that it does not imply any estimate  of $\nu( B(x,r)\cap Y )$ from below in terms of $r^{\upsilon_Y}$ as in a so called lower $\upsilon_Y$-Ahlfors regularity condition  that together with (\ref{defn:uareg1}) would imply the validity
of the so-called $\upsilon_Y$-Ahlfors regularity condition and accordingly the validity of the so-called doubling  condition for the measure $\nu$, \textit{i.e.}, the following condition
\begin{eqnarray} \nonumber
&&\text{there\ exist}\ r_{X,Y}\in]0,+\infty]\,,\ c_{X,Y}\in]0,+\infty[\ \text{such\ that}
\\ \nonumber
&&\nu( B(x,2r)\cap Y )\leq c_{X,Y} \nu( B(x,r)\cap Y )
\\  \label{defn:doubling}
&&\text{for\ all}\ x\in X\ \text{and}\  r \in]0,r_{X,Y}[
 \,. 
\end{eqnarray}
Now condition (\ref{defn:uareg1}) says that  one can estimate the $\nu$-measure of a ball $B(x,r)\cap Y$ in $Y$ from above in terms of the measure of a ball of radius $r$ in a Euclidean space of dimension $\upsilon_Y$ (at least for integer values of $\upsilon_Y$).

As in \cite{La22a}, we also  consider a stronger condition than (\ref{defn:uareg1}) that still does not involve estimates from below for $\nu$  in which we replace the ball $B(x,r)\cap Y$ in $Y$ with an annular domain 
$(B(x,r_2)\setminus B(x,r_1))\cap Y $ with $0\leq r_1< r_2$ in $Y$ and that  says that  one can estimate the $\nu$-measure of an annular domain 
$(B(x,r_2)\setminus B(x,r_1))\cap Y$ in $Y$ from above  in terms of the measure of an annular domain of radii $r_1$ and $r_2$ in a Euclidean space of dimension $\upsilon_Y$ (at least for integer values of $\upsilon_Y$).  Namely, 
  we assume that $Y$ is strongly upper $\upsilon_Y$-Ahlfors regular with respect to $X$, \textit{i.e.}, that
\begin{eqnarray} \nonumber
&&\text{there\ exist}\ r_{X,Y,\upsilon_Y}\in]0,+\infty]\,,\ c_{X,Y,\upsilon_Y}\in]0,+\infty[\ \text{such\ that}
\\ \nonumber
&&\nu( (B(x,r_2)\setminus B(x,r_1))\cap Y )\leq c_{X,Y,\upsilon_Y}(r_2^{\upsilon_Y}-r_1^{\upsilon_Y})
\\ \label{defn:suareg1}
&&\text{for\ all}\ x\in X\ \text{and}\ r_1,r_2\in[0,r_{X,Y,\upsilon_Y}[
\ \text{with}\ r_1<r_2\,,
\end{eqnarray}
where we understand that $B(x,0)\equiv\emptyset$  (in case $X=Y$, we just say that $Y$ is strongly upper $\upsilon_Y$-Ahlfors regular). So, for example,	
 if $Y$ is the boundary of a bounded open Lipschitz  subset of $M={\mathbb{R}}^n$, then $Y$ is upper $(n-1)$-Ahlfors regular with respect to ${\mathbb{R}}^n$ (cf.~Proposition \ref{prop:lgar} of the Appendix \ref{appendix})
and if $Y$ is the boundary of an open bounded subset of $M={\mathbb{R}}^n$ of class $C^1$, then $Y$ is strongly upper $(n-1)$-Ahlfors regular  with respect to $Y$ (cf.~Proposition \ref{prop:lgar} and Remark \ref{rem:coc1sa} of the Appendix \ref{appendix}).

One may wonder about the importance of considering only growth conditions from above for the measure as in the 
   upper Ahlfors or strong upper Ahlfors regularity condition and not from below and about the importance of considering measures that do not satisfy the doubling condition. Here we mention the papers of Verdera \cite{Ve13}, \cite[p.~21]{Ve23} in connection with the integral operator of the Cauchy kernel. Also, in  section \ref{sec:exavio}, we  present some elementary examples of surfaces in ${\mathbb{R}}^3$ in which such conditions from below and the doubling are actually violated.

In the present survey paper $X$ and $Y$ are mainly subsets of the boundary $\partial\Omega$ of some bounded open subset $\Omega$ of $M\equiv {\mathbb{R}}^n$, $d$ is the Euclidean distance and $\nu$ coincides with the ordinary $(n-1)$-dimensional measure on $\partial\Omega$.

\section{Preliminaries and Notation} 
Let $M_n({\mathbb{R}})$ denote the set of $n\times n$ matrices with real entries. $|A|$ denotes the operator norm of a matrix $A$, 
       $A^{t}$ denotes the transpose matrix of $A$.  Let $O_{n}({\mathbb{R}})$ denote the set of $n\times n$ orthogonal matrices with real entries. We also set
\[
 {\mathbb{B}}_{n}(x,\rho)\equiv \left\{
y\in {\mathbb{R}}^{n}:\,\vert x-y\vert <\rho
\right\}   
\qquad \forall
(\xi,\rho)\in {\mathbb{R}}^n\times ]0,+\infty[
\,.
\]
 Here and in the sequel, $m_n$ denotes the $n$-dimensional Lebesgue measure in ${\mathbb{R}}^n$ and  $m_{n-1}$ denotes the ordinary $(n-1)$-dimensional (surface) measure and
\begin{equation}\label{eq:snon}
\omega_n\equiv m_n({\mathbb{B}}_n(0,1))\,,\qquad
s_n\equiv m_{n-1}(\partial{\mathbb{B}}_n(0,1))\,.
\end{equation}
For the standard notation of the spaces of H\"{o}lder or Lipschitz continuous functions, we refer for example  to \cite[\S 2]{DoLa17}, \cite[\S 2.6]{DaLaMu21}.  Let $\Omega$ be an open subset of ${\mathbb{R}}^n$.  Let $s\in {\mathbb{N}}\setminus\{0\}$, $f\in \left(C^{1}(\Omega)\right)^{s} $. Then   $Df$ denotes the Jacobian matrix of $f$. In order to analyze the kernel of the double layer potential, we need some more information on the fundamental solution $S_{ {\mathbf{a}} } $. To do so, we introduce the fundamental solution $S_{n}$ of the Laplace operator. Namely, we set
\[
%\label{ups}
S_{n}(x)\equiv
\left\{
\begin{array}{lll}
\frac{1}{s_{n}}\ln  |x| \qquad &   \forall x\in 
{\mathbb{R}}^{n}\setminus\{0\},\quad & {\mathrm{if}}\ n=2\,,
\\
\frac{1}{(2-n)s_{n}}|x|^{2-n}\qquad &   \forall x\in 
{\mathbb{R}}^{n}\setminus\{0\},\quad & {\mathrm{if}}\ n>2\,.
\end{array}
\right.
\]
 and
we follow a formulation of Dalla Riva \cite[Thm.~5.2, 5.3]{Da13} and Dalla Riva, Morais and Musolino \cite[Thm.~5.5]{DaMoMu13}, that we state  as in paper  \cite[Cor.~4.2]{DoLa17} of Dondi and the author  (see also John~\cite{Jo55}, and Miranda~\cite{Mi65} for homogeneous operators, and Mitrea and Mitrea~\cite[p.~203]{MitMit13}).   
 \begin{proposition}
 \label{prop:ourfs} 
Let ${\mathbf{a}}$ be as in (\ref{introd0}), (\ref{ellip}), (\ref{symr}). 
Let $S_{ {\mathbf{a}} }$ be a fundamental solution of $P[{\mathbf{a}},D]$. 
Then there exist an invertible matrix $T\in M_{n}({\mathbb{R}})$ such that
\begin{equation}
\label{prop:ourfs0}
a^{(2)}=TT^{t}\,,
\end{equation}
 a real analytic function $A_{1}$ from $\partial{\mathbb{B}}_{n}(0,1)\times{\mathbb{R}}$ to ${\mathbb{C}}$ such that
 $A_{1}(\cdot,0)$ is odd,    $b_{0}\in {\mathbb{C}}$, a real analytic function $B_{1}$ from ${\mathbb{R}}^{n}$ to ${\mathbb{C}}$ such that $B_{1}(0)=0$, and a real analytic function $C $ from ${\mathbb{R}}^{n}$ to ${\mathbb{C}}$ such that
\begin{eqnarray}%attenzione a formulare diversamente nel lavoro
\label{prop:ourfs1}
\lefteqn{
S_{ {\mathbf{a}} }(x)
= 
\frac{1}{\sqrt{\det a^{(2)} }}S_{n}(T^{-1}x)
}
\\ \nonumber
&&\qquad
+|x|^{3-n}A_{1}(\frac{x}{|x|},|x|)
 +(B_{1}(x)+b_{0}(1-\delta_{2,n}))\ln  |x|+C(x)\,,
 %\\ \nonumber
%&&\qquad
%-\frac{\delta_{2,n}}{2\pi}\int_{\partial{\mathbb{B}}_{2}(0,1)}S_{a^{(2)}}-A_{0}\,d\sigma
\end{eqnarray}
for all $x\in {\mathbb{R}}^{n}\setminus\{0\}$,
 and such that both $b_{0}$ and $B_{1}$   equal zero
if $n$ is odd. Moreover, 
 \[
 \frac{1}{\sqrt{\det a^{(2)} }}S_{n}(T^{-1}x) 
 \]
is a fundamental solution for the principal part
  of $P[{\mathbf{a}},D]$.
\end{proposition}
In particular for the statement that $A_{1}(\cdot,0)$ is odd, we refer to
Dalla Riva, Morais and Musolino \cite[Thm.~5.5, (32)]{DaMoMu13}, where $A_{1}(\cdot,0)$ coincides with ${\mathbf{f}}_1({\mathbf{a}},\cdot)$ in that paper.  Then for each $\theta\in]0,1]$, we define the function $\omega_{\theta}(\cdot)$ from $[0,+\infty[$ to itself by setting
\begin{equation}\label{omth}
\omega_{\theta}(r)\equiv
\left\{
\begin{array}{ll}
0 &r=0\,,
\\
r^{\theta}\vert \ln r \vert  &r\in]0,r_{\theta}]\,,
\\
r_{\theta}^{\theta}\vert \ln r_{\theta} \vert  & r\in ]r_{\theta},+\infty[\,,
\end{array}
\right.
\end{equation}
where
$
%\label{omth1}
r_{\theta}\equiv e^{-1/\theta}
$. 
Next we introduce some notation for the kernels. We do so in the abstract context of metric spaces. If $X$ and $Y$ are subsets of a metric space $M$, we consider off-diagonal kernels $K$ from $(X\times Y)\setminus D_{X\times Y}$ to ${\mathbb{C}}$, where
\[
D_{X\times Y}\equiv  \left\{
(x,y)\in X\times Y:\,x=y
\right\}
\] denotes the diagonal set of $X\times Y$  
 and we introduce the following class of   `potential type' kernels  (see also paper \cite{DoLa17} of the author and Dondi,  where such classes have been introduced in a form that generalizes those of Giraud \cite{Gi34}, Gegelia \cite{Ge67}, 
 Kupradze, Gegelia, Basheleishvili and 
 Burchuladze \cite[Chap.~IV]{KuGeBaBu79}).	 
\begin{definition}\label{defn:ksss} Let 
$
(M,d)$ be a metric space.
 Let $X$, $Y\subseteq M$. Let $s_1$, $s_2$, $s_3\in {\mathbb{R}}$. We denote by the symbol ${\mathcal{K}}_{s_1, s_2, s_3} (X\times Y)$ the set of continuous functions $K$ from $(X\times Y)\setminus D_{X\times Y}$ to ${\mathbb{C}}$ such that
 \begin{eqnarray*}
\lefteqn{
\|K\|_{  {\mathcal{K}}_{ s_1, s_2, s_3  }(X\times Y)  }
\equiv
\sup\biggl\{\biggr.
d(x,y)^{ s_{1} }\vert K(x,y)\vert :\,(x,y)\in X\times Y, x\neq y
\biggl.\biggr\}
}
\\ \nonumber
&&\qquad\qquad\qquad
+\sup\biggl\{\biggr.
\frac{d(x',y)^{s_{2}}}{d(x',x'')^{s_{3}}}
\vert   K(x',y)- K(x'',y)  \vert :\,
\\ \nonumber
&&\qquad\qquad\qquad 
x',x''\in X, x'\neq x'', y\in Y\setminus B(x',2d(x',x''))
\biggl.\biggr\}<+\infty\,.
\end{eqnarray*}
\end{definition}
For  $s_2=s_1+s_3$ one has the so-called class of standard kernels  that is the case in which Garc\'{\i}a-Cuerva and Gatto \cite{GaGa04}, \cite{GaGa05}, Gatto \cite{Gat09} have proved  $T1$ Theorems for the integral operators with kernel $K$ in case of weakly singular, singular and hyper-singular integral operators with $X=Y$.\par 

We also consider the following more restrictive class of kernels.
\begin{definition}Let 
$
(M,d)$ be a metric space. Let $X$, $Y\subseteq M$. Let $\nu$ be as in (\ref{eq:nu}).   Let $s_1$, $s_2$, $s_3\in {\mathbb{R}}$. We set
\begin{eqnarray*}\nonumber
\lefteqn{
 {\mathcal{K}}_{ s_1, s_2, s_3  }^\sharp(X\times Y)
  \equiv 
  \biggl\{\biggr.
K\in  {\mathcal{K}}_{ s_1, s_2, s_3  }(X\times Y):\,
}
\\ \nonumber
&&\ \ 
K(x,\cdot)\ \text{is}\ \nu-\text{integrable\ in}\ Y\setminus B(x,r)
\ \text{for\ all}\ (x,r)\in X\times]0,+\infty[\,,
\\ \nonumber
&&\ \ 
\sup_{x\in X}\sup_{r\in ]0,+\infty[}
\left\vert 
\int_{Y\setminus B(x,r)}K(x,y)\,d\nu(y)\right\vert<+\infty
 \biggl.\biggr\}  
\end{eqnarray*}
and
\begin{eqnarray*}
\lefteqn{
\|K\|_{{\mathcal{K}}_{ s_1, s_2, s_3  }^\sharp(X\times Y)}
\equiv
\|K\|_{{\mathcal{K}}_{ s_1, s_2, s_3  }(X\times Y)}
}
\\ \nonumber
&&+
\sup_{x\in X}\sup_{r\in ]0,+\infty[}
\left\vert
\int_{Y\setminus B(x,r)}K(x,y)\,d\nu(y)
\right\vert\quad\forall 
 K\in {\mathcal{K}}_{ s_1, s_2, s_3  }^\sharp(X\times Y)\,.
 \end{eqnarray*}
\end{definition}
Clearly,  $({\mathcal{K}}^\sharp_{ s_{1},s_{2},s_{3}   }(X\times Y),\|\cdot\|_{  {\mathcal{K}}^\sharp_{s_{1},s_{2},s_{3}   }(X\times Y)  })$ is a normed space and the space ${\mathcal{K}}^\sharp_{ s_{1},s_{2},s_{3}   }(X\times Y)$ is continuously embedded into ${\mathcal{K}}_{ s_{1},s_{2},s_{3}   }(X\times Y)$.\par

\section{The double layer potential on the boundary of Lipschitz and $C^1$ domains}
\label{dolalipc1}

For the definition of   open Lipschitz subset of ${\mathbb{R}}^n$ and of open subset of ${\mathbb{R}}^n$ of class $C^1$, we refer for example to 
Dalla Riva, the author and Musolino \cite[\S 2.9, 2.13]{DaLaMu21}.

Salaev \cite{Sa76} has proved that the Cauchy integral on a rectifiable simple closed curve that satisfies an upper Ahlfors regularity condition  is bounded in the space  $C^{0,\beta}$ on the curve for $\beta\in]0,1[$.
 
 Since the double layer potential is the real part of the Cauchy integral for real densities, it follows that the double layer potential  $W_\Omega[{\mathbf{a}},S_{{\mathbf{a}}},\cdot]$ is bounded in $C^{0,\beta}(\partial\Omega)$ for $\beta\in]0,1[$  in case $\Omega$ is a Jordan domain
 bounded by a rectifiable simple closed curve that satisfies an upper Ahlfors regularity condition and $S_{{\mathbf{a}}}$ equals the fundamental solution of the Laplace operator.
 
 For the estimate of moduli of continuity  of the Cauchy integral on a rectifiable simple closed curve that satisfies an upper Ahlfors regularity condition even in generalized H\"{o}lder spaces, we should also mention the papers of   Plemelj
\cite{Plemelj-6-2},  Privaloff \cite{Priv16-6-2},  Zygmund
\cite{Zyg02-6-2},  Magnaradze \cite{Magn-6-2},  Babaev and
 Salaev \cite{Bab-Sal73-6-2},  Tamrazov
\cite{Tam-mono-6-2}, \cite{Tam78a-6-2}, \cite{Tam78-6-2},    Gerus
\cite{G77-6-2}, \cite{G78-6-2}, \cite{G96-6-2}, \cite{Gerus-6-2},  Salimov
\cite{Salimov-6-2},  Dyn'kin \cite{Dynkin-6-2},  Salaev,
 Guse\u{\i}nov and  Se\u{\i}fullaev \cite{SaGusSe90},
 Guse\u{\i}nov \cite{Gus92}.

 Then Mitrea,  Mitrea and Mitrea \cite[Prop.~25.5.21]{MitMitMit23a} have proved that the double layer potential  $W_\Omega[{\mathbf{a}},S_{{\mathbf{a}}},\cdot]$ is bounded in $C^{0,\beta}(\partial\Omega)$  for $\beta\in]0,1[$  in case $S_{{\mathbf{a}}}$ equals the fundamental solution of the Laplace operator and 
the ordinary $(n-1)$-dimensional measure on the boundary on  $\Omega$ satisfies an upper Ahlfors growth condition, a condition that certainly holds if $\Omega$ is
 a  bounded open  Lipschitz subset of ${\mathbb{R}}^n$ with $n\geq 2$.  
 
 If $\Omega$ is a bounded open Lipschitz subset of ${\mathbb{R}}^n$, then there exists a   subset $N$  of measure zero of $\partial\Omega$ such that the outward unit normal $\nu$ exists at all points of $(\partial\Omega)\setminus N$.

By Mitrea,  Mitrea and Mitrea \cite[Prop.~25.5.21]{MitMitMit23a}, Mitrea \cite{Mit23}, the function $W_\Omega[{\mathbf{a}},S_{{\mathbf{a}}}   ,\mu]$ is uniformly continuous on $ (\partial\Omega)\setminus N$ and thus it admits a unique uniformly continuous extension
 to the closure of $(\partial\Omega)\setminus N$ in $\partial\Omega$, \textit{i.e.}, to the whole of $\partial\Omega$. However, the principal value of (\ref{introd3}) may well exist at some point $x\in N$ and be different from the value of the continuous extension of $W_\Omega[{\mathbf{a}},S_{{\mathbf{a}}}   ,\mu]$ at $x$. For the existence of the principal value of (\ref{introd3}) at points of $N$ in dimension 2, we refer to classical books such as that of Gakhov \cite[\S 4.5, p.~31]{Ga66} and in higher dimensions we refer to     Burago and Maz'ya~\cite[Thm.~2, p.~17]{BuMa69}. In case $\Omega$ is of class $C^1$, we can take $N=\emptyset$.

 \section{The double layer potential  on the boundary of domains of class $C^{1,\alpha}$ with $\alpha\in]0,1[$.}
 
 For the definition of    open subset of ${\mathbb{R}}^n$ of class $C^{1,\alpha}$ for some $\alpha\in]0,1]$, we refer for example to 
Dalla Riva, the author and Musolino \cite[\S 2.13]{DaLaMu21}.  
 
 Here we must say that we only consider the boundary behaviour of the double layer potential. Instead for the regularity properties of the double layer potential in the  Schauder space $C^{1,\alpha}$ outside of the boundary we refer to G\"{u}nter~\cite{Gu67}, Kupradze,  Gegelia,  Basheleishvili and 
 Burchuladze~\cite{KuGeBaBu79}, Mikhlin \cite{Mik70}, 
Mikhlin and  Pr\"{o}ssdorf \cite{MikPr86}, Miranda \cite{Mi65},   \cite{Mi70},  
Wiegner~\cite{Wi93}, Dalla Riva \cite{Da13}, Dalla Riva, Morais and Musolino \cite{DaMoMu13}, Mitrea, Mitrea and Verdera \cite{MitMitVe16} and references therein.

 We first state the following theorem, that  extends a known result of Schauder \cite[p.~614]{Sc31} for the harmonic double layer potential in case $n=3$. For later contributions see also Mitrea   \cite{Mit14}. For the (classical) definition of the generalized H\"{o}lder space    $C^{0,\omega_1(\cdot)}(\partial\Omega)$ we refer for example to \cite[\S 2]{DoLa17}.
 \begin{theorem}\label{thm:dllregen}
 Let $n\in {\mathbb{N}}\setminus\{0,1\}$. 
 Let ${\mathbf{a}}$ be as in (\ref{introd0}), (\ref{ellip}), (\ref{symr}).  Let $S_{ {\mathbf{a}} }$ be a fundamental solution of $P[{\mathbf{a}},D]$.   Let $\alpha\in]0,1[$, $\beta\in]0,1]$. 
 
 %In case $ (n,\alpha)=(2,1)$, we further assume that $DB_1(0)=0$.
  
 Let $\Omega$ be a bounded open subset of ${\mathbb{R}}^{n}$ of class $C^{1,\alpha}$.   Then the following statements hold.
 \begin{enumerate}
 \item[(i)] If $0<\beta<1-\alpha $, then the    operator $W_\Omega[{\mathbf{a}},S_{{\mathbf{a}}}   ,\cdot]$ from 
$C^{0,\beta}(\partial\Omega)$ to $C^{0,\alpha+\beta}(\partial\Omega)$ defined by (\ref{introd3})
 for all $\mu\in C^{0,\beta}(\partial\Omega)$ is linear and continuous.
 \item[(ii)] If $ \beta=1-\alpha $, then the    operator $W_\Omega[{\mathbf{a}},S_{{\mathbf{a}}}   ,\cdot]$ from 
$C^{0,\beta}(\partial\Omega)$ to $C^{0,\omega_1(\cdot)}(\partial\Omega)$ defined by (\ref{introd3})
 for all $\mu\in C^{0,\beta}(\partial\Omega)$ is linear and continuous.

  \end{enumerate}
\end{theorem} 
{\bf Proof.}  We first note that the membership of $1$ in $C^{1,\alpha}(\partial\Omega)$ and Dondi and the author  \cite[Thm 9.2]{DoLa17}  implies that 
\[
W_\Omega[{\mathbf{a}},S_{{\mathbf{a}}}   ,1]\in C^{1,\alpha}(\partial\Omega)\subseteq C^{0,\omega_1(\cdot)}(\partial\Omega)\subseteq C^{0,\theta}(\partial\Omega)\qquad\forall\theta\in]0,1[\,.
\]
 By \cite[Rmk 6.1 (ii)]{DoLa17}, we know that the kernel $\overline{B^{*}_{\Omega,y}}\left(S_{{\mathbf{a}}}(x-y)\right)$ of the double layer potential belongs to the class
${\mathcal{K}}_{n-1-\alpha,n-\alpha,1}((\partial\Omega)\times(\partial\Omega))$. In particular, $\overline{B^{*}_{\Omega,y}}\left(S_{{\mathbf{a}}}(x-y)\right)$ is weakly singular and accordingly $W_\Omega[{\mathbf{a}},S_{{\mathbf{a}}}   ,\cdot]$ defines a linear and continuous map from $C^{0,\beta}(\partial\Omega)$ to $C^{0}(\partial\Omega)$ both in case of statement (i) and of statement (ii) (cf.~\textit{e.g.}, 
\cite[Prop.~6.1 (i)]{DoLa17}).
We also note that in case $(n-\alpha)-\beta>n-1$ of statement (i), we also have
$1+(n-1)-(n-\alpha-\beta)>0$. Then \cite[Prop.~5.11]{La22a} (or 
\cite[Lem.~6.1]{DoLa17}) implies the existence of $c\in]0,+\infty[$ such that
\begin{eqnarray}
\label{k0b1}
\lefteqn{
| W_\Omega[{\mathbf{a}},S_{{\mathbf{a}}}   ,\mu](x')-W_\Omega[{\mathbf{a}},S_{{\mathbf{a}}}   ,\mu](x'')|
}
\\ \nonumber
&&\qquad
\leq
c\| K \|_{  {\mathcal{K}}_{n-1-\alpha,n-\alpha,1}((\partial\Omega)\times(\partial\Omega))  } \|\mu\|_{  C^{0,\beta}(\partial\Omega)  }\omega(|x'-x''|)
\\ \nonumber
&&\qquad\quad
+\|\mu\|_{  C^{0}(\partial\Omega)  }
| W_\Omega[{\mathbf{a}},S_{{\mathbf{a}}}   ,1](x')-W_\Omega[{\mathbf{a}},S_{{\mathbf{a}}}   ,1](x'')|
\end{eqnarray}
for all $x',x''\in\partial\Omega$ such that $|x'-x''|<e^{-1}$ and for all  $ \mu\in C^{0,\beta}(\partial\Omega)$, where 
\[
\omega(r)\equiv \left\{
\begin{array}{ll}
r^{\min\{\alpha+\beta,1\}}  & {\mathrm{if}}\ n-\alpha-\beta<n-1\ \text{as\ in}\ (i)\,,
\\
\max\{
r^{\alpha+\beta},\omega_1(r)\}  & {\mathrm{if}}\ n-\alpha-\beta=n-1\ \text{as\ in}\ (ii)\,,
\end{array}
\right.
\ \ \forall r\in]0,+\infty[\,.
\]
Then under the assumptions of statement (i), we have
\[
C^{0,r^{\min\{\alpha+\beta,1\}} }(\partial\Omega)=C^{0,\alpha+\beta}(\partial\Omega)
\]
and the membership of $W_\Omega[{\mathbf{a}},S_{{\mathbf{a}}}   ,1]$ in $C^{0,\alpha+\beta}(\partial\Omega)$ and inequality (\ref{k0b1}) imply the validity of statement (i) (see also \cite[Rmk. 2.2]{DoLa17}).  Instead, under the assumptions of statement (ii), we have $\alpha+\beta=1$, 
\[
C^{0,\max\{
r^{\alpha+\beta},\omega_1(r)\} }(\partial\Omega)=C^{0,\omega_1(\cdot)}(\partial\Omega)
\]
and the membership of $W_\Omega[{\mathbf{a}},S_{{\mathbf{a}}}   ,1]$ in $C^{0,\omega_1(\cdot)}(\partial\Omega)$ and inequality (\ref{k0b1}) imply the validity of statement (ii) (see also \cite[Rmk. 2.2]{DoLa17}).\hfill  $\Box$ 

\vspace{\baselineskip}

Next we turn to consider case $\alpha+\beta>1$ and we state the following theorem, that collects and extends results of 
 Fichera and De Vito \cite[LXXXIII]{FiDe70} for the Laplace operator in case $n=2$. See also  Miranda \cite[15.VI]{Mi70}, where the author mentions a result of Giraud \cite{Gi32}, Mitrea~\cite{Mit14},  Dondi and the author  \cite{DoLa17} and   \cite[Thm.~5.1]{La22d}.    
 For the (classical) definition of the generalized Schauder space    $C^{1,\omega_\alpha(\cdot)}(\partial\Omega)$ we refer for example to \cite[\S 2]{DoLa17}.
 \begin{theorem}\label{thm:dllregenn}
 Let $n\in {\mathbb{N}}\setminus\{0,1\}$. 
 Let ${\mathbf{a}}$ be as in (\ref{introd0}), (\ref{ellip}), (\ref{symr}).  Let $S_{ {\mathbf{a}} }$ be a fundamental solution of $P[{\mathbf{a}},D]$.   Let $\alpha\in]0,1[$, $\beta\in]0,1]$, $\alpha+\beta>1$. 
 
 %In case $ (n,\alpha)=(2,1)$, we further assume that $DB_1(0)=0$.
  
 Let $\Omega$ be a bounded open subset of ${\mathbb{R}}^{n}$ of class $C^{1,\alpha}$.   Then the following statements hold.
 \begin{enumerate}
\item[(i)] If $  \beta<1$, then the    operator $W_\Omega[{\mathbf{a}},S_{{\mathbf{a}}}   ,\cdot]$ from 
$C^{0,\beta}(\partial\Omega)$ to $C^{1,\alpha+\beta-1}(\partial\Omega)$ defined by (\ref{introd3})
 for all $\mu\in C^{0,\beta}(\partial\Omega)$ is linear and continuous.
\item[(ii)] If $\beta=1$, then the  operator $W_\Omega[{\mathbf{a}},S_{{\mathbf{a}}}   ,\cdot]$
 from 
$C^{0,\beta}(\partial\Omega)=C^{0,1}(\partial\Omega)$ to $C^{1,\omega_{\alpha+\beta-1}}(\partial\Omega)=C^{1,\omega_\alpha}(\partial\Omega)$ defined by (\ref{introd3})
 for all $\mu\in  C^{0,1}(\partial\Omega)$ 
 is linear and continuous.
 \end{enumerate}
\end{theorem}
For a proof we refer to \cite[Thm.~5.1]{La22d}. Here we do not provide a complete proof of Theorem \ref{thm:dllregenn}, but we  point out that the proof is based on statements that hold in the general setting that we have mentioned in the introduction. Indeed, $\partial\Omega$ is   strongly upper $(n-1)$-Ahlfors regular (cf.~Proposition \ref{prop:lgar} and Remark \ref{rem:coc1sa} of the Appendix \ref{appendix}) and Dondi and the author  \cite[Thm 9.2]{DoLa17}  implies that 
\begin{equation}\label{eq:dllregenn}
W_\Omega[{\mathbf{a}},S_{{\mathbf{a}}}   ,1]\in C^{1,\alpha}(\partial\Omega) 
\end{equation}
and we note that 
  one can prove classically that
\begin{eqnarray}\label{eq:gradformn}
\lefteqn{
{\mathrm{grad}}_{\partial\Omega,x}W_\Omega[{\mathbf{a}},S_{{\mathbf{a}}}   ,\mu](x)
=
{\mathrm{grad}}_{\partial\Omega,x} \int_{\partial\Omega}  
\overline{B^{*}_{\Omega,y}}\left(S_{{\mathbf{a}}}(x-y)\right)\mu(y)\,d\sigma_y
}
\\ \nonumber
&&\qquad\qquad\qquad\qquad
=\int_{\partial\Omega}[{\mathrm{grad}}_{\partial\Omega,x} \overline{B^{*}_{\Omega,y}}\left(S_{{\mathbf{a}}}(x-y)\right)]
 (\mu(y)-\mu(x))\,d\sigma_y
\\ \nonumber
&&\qquad\qquad\qquad\qquad\quad
+\mu(x){\mathrm{grad}}_{\partial\Omega} \int_{\partial\Omega} 
 \overline{B^{*}_{\Omega,y}}\left(S_{{\mathbf{a}}}(x-y)\right)
 \,d\sigma_y 
 \\ \nonumber
&&\qquad\qquad\qquad\qquad
=
 Q[K,\mu,1](x)
 \\ \nonumber
&&\qquad\qquad\qquad\qquad\quad
 +\mu(x){\mathrm{grad}}_{\partial\Omega}W_\Omega[{\mathbf{a}},S_{{\mathbf{a}}}   ,1](x)
 \qquad\forall x\in\partial\Omega 
 \end{eqnarray}
for all $\mu\in C^{0,\beta}(\partial\Omega)$, where ${\mathrm{grad}}_{\partial\Omega}$ denotes the tangental gradient and
${\mathrm{grad}}_{\partial\Omega,x}$ denotes the tangental gradient with respect to the first variable (cf.~\cite[Thm.~6.1]{La22b}) and
\begin{eqnarray*}
&&
K(x,y)\equiv
-[{\mathrm{grad}}_{\partial\Omega,x} \overline{B^{*}_{\Omega,y}}\left(S_{{\mathbf{a}}}(x-y)\right)]\quad\forall (x,y)\in (\partial\Omega)^2\setminus D_{(\partial\Omega)\times (\partial\Omega)}\,,
\\ \nonumber
&&
Q[K,\mu,1](x)
\equiv 
\int_{\partial\Omega} K(x,y)
 (\mu(x)-\mu(y))\,d\sigma_y \quad\forall x\in\partial\Omega 
\end{eqnarray*}
for all $\mu\in C^{0,\beta}(\partial\Omega)$.
 Then one can deduce  the continuity of the operator $Q[K,\cdot,1]$  from $C^{0,\beta}(\partial\Omega)$ to $C^{0,\alpha+\beta-1}(\partial\Omega)$ if $\beta<1$ and to $C^{0,\omega_\alpha}(\partial\Omega)$ if $\beta=1$ by exploiting the proof of  \cite[Lem.~4.6]{La22d} on the membership of the tangential gradient of the kernel of the double layer potential in an appropriate class and 
the abstract result on $Q$ of \cite[Prop.~6.3 (iii) (c), (cc)]{La22a} in metric spaces, that generalizes previous work of Gatto \cite{Gat09}. Then  the membership of (\ref{eq:dllregenn})
 together with the  continuity of the pointwise product in generalized Schauder spaces (cf.~\textit{e.g.}, \cite[Lems.~2.4, 2.5]{DoLa17}) imply that ${\mathrm{grad}}_{\partial\Omega,x}W_\Omega[{\mathbf{a}},S_{{\mathbf{a}}}   ,\cdot]$ is bounded from $C^{0,\beta}(\partial\Omega)$ to $C^{1,\alpha+\beta-1}(\partial\Omega)$ if $\beta<1$ and to $C^{1,\omega_\alpha}(\partial\Omega)$ if $\beta=1$ and thus   Theorem \ref{thm:dllregenn}  can be proved to be true.\par

By setting $\beta=\alpha$ in the previous Theorems \ref{thm:dllregen}, \ref{thm:dllregenn}, we immediately deduce the validity of the following corollary that says that in a set of a class $C^{1,\alpha}$ with $\alpha\in]0,1[$. the regularizing effect of  boundary double layer potential on the boundary equals $\alpha$, with the only exceptional value $\alpha=1/2$.
\begin{corollary}
 Let $n\in {\mathbb{N}}\setminus\{0,1\}$. 
 Let ${\mathbf{a}}$ be as in (\ref{introd0}), (\ref{ellip}), (\ref{symr}).  Let $S_{ {\mathbf{a}} }$ be a fundamental solution of $P[{\mathbf{a}},D]$.   Let $\alpha\in]0,1[$. 
 
 %In case $ (n,\alpha)=(2,1)$, we further assume that $DB_1(0)=0$.
  
 Let $\Omega$ be a bounded open subset of ${\mathbb{R}}^{n}$ of class $C^{1,\alpha}$.   Then the following statements hold.
 \begin{enumerate}
 \item[(i)] If $0<\alpha<1/2 $, then the    operator $W_\Omega[{\mathbf{a}},S_{{\mathbf{a}}}   ,\cdot]$ from 
$C^{0,\alpha}(\partial\Omega)$ to $C^{0,2\alpha}(\partial\Omega)$ defined by (\ref{introd3})
 for all $\mu\in C^{0,\alpha}(\partial\Omega)$ is linear and continuous.
 \item[(ii)] If $\alpha=1/2 $, then the    operator $W_\Omega[{\mathbf{a}},S_{{\mathbf{a}}}   ,\cdot]$ from 
$C^{0,\alpha}(\partial\Omega)$ to $C^{0,\omega_{2\alpha}(\cdot)}(\partial\Omega)$ defined by (\ref{introd3})
 for all $\mu\in C^{0,\alpha}(\partial\Omega)$ is linear and continuous.
 \item[(iii)] If $1/2<\alpha<1$, then the    operator $W_\Omega[{\mathbf{a}},S_{{\mathbf{a}}}   ,\cdot]$ from 
$C^{0,\alpha}(\partial\Omega)$ to $C^{1,2\alpha-1}(\partial\Omega)$ defined by (\ref{introd3})
 for all $\mu\in C^{0,\alpha}(\partial\Omega)$ is linear and continuous.
\end{enumerate}
\end{corollary}

 \section{The double layer potential   on the boundary of domains of class $C^{1,1}$ and $C^2$.}
 
 We first mention some known results in the classical case of the boundary behaviour of the  double layer potential in Schauder spaces with $m= 2$. Instead for the regularity properties of the double layer potential in Schauder spaces with $m= 2$ outside of the boundary we refer to G\"{u}nter~\cite{Gu67}, Kupradze,  Gegelia,  Basheleishvili and 
 Burchuladze~\cite{KuGeBaBu79}, Mikhlin \cite{Mik70}, 
Mikhlin and  Pr\"{o}ssdorf \cite{MikPr86}, Miranda \cite{Mi65},   \cite{Mi70},  
Wiegner~\cite{Wi93}, Dalla Riva \cite{Da13}, Dalla Riva, Morais and Musolino \cite{DaMoMu13}, Mitrea, Mitrea and Verdera \cite{MitMitVe16} and references therein.\par
 
  In case $n=3$  and $\Omega$ is of class $C^{2}$, $\alpha\in]0,1[$   and if  $P[{\mathbf{a}},D]$ is the  Helmholtz operator, Colton and Kress~\cite{CoKr83} have 
developed previous work of G\"{u}nter~\cite{Gu67} and Mikhlin~\cite{Mik70}  and proved that the operator $W[\partial\Omega ,{\mathbf{a}},S_{{\mathbf{a}}},\cdot]$
is bounded from $C^{0,\alpha}(\partial\Omega)$ to $C^{1,\alpha}(\partial\Omega)$.

In case $n\geq 2$, $\alpha\in]0,1[$  and $\Omega$ is of class $C^{2}$   and if  $P[{\mathbf{a}},D]$ is the  Laplace operator, 
Hsiao and Wendland \cite[Remark 1.2.1, p.~10]{HsWe08} 
deduce  that the operator $W[\partial\Omega ,{\mathbf{a}},S_{{\mathbf{a}}},\cdot]$
is bounded from $C^{0,\alpha}(\partial\Omega)$ to $C^{1,\alpha}(\partial\Omega)$ by the work of
Mikhlin and Pr\"{o}ssdorf \cite{MikPr86}.

 We now show that if $\Omega$ is of class  $C^{1,1}$, then the double layer potential improves the regularity of one unit if $\beta<1$ (and of less than that if $\beta=1$) Namely, we have the following statement that  is a generalization of   classical results  for the Laplace and Helmoltz operator in case $\Omega$ is of class $C^2$ and $\beta\in]0,1[$ (see Colton and Kress \cite[Thm. 2.22]{CoKr83}, 
Hsiao and Wendland \cite[Remark 1.2.1, p.~10]{HsWe08}).   \begin{theorem}\label{thm:dllreggen}
  Let  $\beta\in]0,1]$.  Let $\Omega$ be a bounded open subset of ${\mathbb{R}}^{n}$ of class  $C^{1,1}$.  
    Let ${\mathbf{a}}$ be as in (\ref{introd0}), (\ref{ellip}), (\ref{symr}).  Let $S_{ {\mathbf{a}} }$ be a fundamental solution of $P[{\mathbf{a}},D]$. Then the following statements hold.
  \begin{enumerate}
\item[(i)] If $\beta<1$, then the  operator $W_\Omega[{\mathbf{a}},S_{{\mathbf{a}}}   ,\cdot]$ from 
$C^{0,\beta}(\partial\Omega)$ to $C^{1, \beta }(\partial\Omega)$ that is defined by 
 (\ref{introd3}) for all $\mu\in C^{0,\beta}(\partial\Omega)$ is linear and continuous.
 \item[(ii)] If $\beta=1$, then the  operator $W_\Omega[{\mathbf{a}},S_{{\mathbf{a}}}   ,\cdot]$ from 
$C^{0,1}(\partial\Omega)$ to $C^{1, \omega_1(\cdot) }(\partial\Omega)$ that is  defined by 
 (\ref{introd3})  for all $\mu\in C^{0,1}(\partial\Omega)$ is linear and continuous.
\end{enumerate}
\end{theorem}
For a proof we refer to \cite[Thm.~1.1]{La23b}.  Here we do not provide a complete proof of Theorem \ref{thm:dllreggen}, but we  point out that the proof is based on statements that hold in the general setting that we have mentioned in the introduction. Indeed, $\partial\Omega$ is   strongly upper $(n-1)$-Ahlfors regular (cf.~Proposition \ref{prop:lgar} and Remark \ref{rem:coc1sa} of the Appendix \ref{appendix}) and reference \cite[Lem. 5.4]{La22d} implies that 
\begin{equation}\label{eq:dllreggen}
W_\Omega[{\mathbf{a}},S_{{\mathbf{a}}}   ,1]\in C^{1,\omega_1(\cdot)}(\partial\Omega)\,.
\end{equation}
and we note that  one can exploit formula (\ref{eq:gradformn}), the proof of  \cite[Lem.~4.6]{La22d} on the membership of the tangential gradient of the kernel of the double layer potential in an appropriate class  and the abstract result of \cite[Prop.~6.3 (ii) (b), (bb)]{La22a} in metric spaces, that generalizes previous work of Gatto \cite{Gat09}. To do so, we deduce the continuity of $Q[K,\cdot,1]$ from $C^{0,\beta}(\partial\Omega)$ to $C^{0, \beta}(\partial\Omega)$ if $\beta<1$ and to $C^{0,\omega_1(\cdot)}(\partial\Omega)$ if $\beta=1$ by the abstract result on $Q$ of \cite[Prop.~6.3 (ii) (b), (bb)]{La22a}, that generalizes previous work of Gatto \cite{Gat09}.  Indeed, $\partial\Omega$ is strongly upper $(n-1)$-Ahlfors regular (cf.~Proposition \ref{prop:lgar} and Remark \ref{rem:coc1sa} of the Appendix \ref{appendix}). Then  the membership of (\ref{eq:dllreggen})  together with the  continuity of the pointwise product in generalized Schauder spaces (cf.~\textit{e.g.}, \cite[Lems.~2.4, 2.5]{DoLa17})  imply that ${\mathrm{grad}}_{\partial\Omega,x}W_\Omega[{\mathbf{a}},S_{{\mathbf{a}}}   ,\cdot]$ is bounded from $C^{0,\beta}(\partial\Omega)$ to $C^{1,\beta}(\partial\Omega)$ if $\beta<1$ and to $C^{1,\omega_1(\cdot)}(\partial\Omega)$ if $\beta=1$ and one can deduce the validity of Theorem \ref{thm:dllreggen}.\par

  \section{The double layer potential   on the boundary of domains of class $C^{m,\alpha}$ with $m\geq 2$, $\alpha\in]0,1[$.}
  We first mention some known results in the classical case of the boundary behaviour of the double layer potential in Schauder spaces with $m\geq 2$. Instead for the regularity properties of the double layer potential in Schauder spaces with $m\geq 2$ outside of the boundary we refer to G\"{u}nter~\cite{Gu67}, Kupradze,  Gegelia,  Basheleishvili and 
 Burchuladze~\cite{KuGeBaBu79}, Mikhlin \cite{Mik70}, 
Mikhlin and  Pr\"{o}ssdorf \cite{MikPr86}, Miranda \cite{Mi65},   \cite{Mi70},  
Wiegner~\cite{Wi93}, Dalla Riva \cite{Da13}, Dalla Riva, Morais and Musolino \cite{DaMoMu13}, Mitrea, Mitrea and Verdera \cite{MitMitVe16} and references therein.\par

  In case  $n=3$, $m\geq 2$, $\alpha\in]0,1]$ and $\Omega$ is of class $C^{m,\alpha}$ and if  $P[{\mathbf{a}},D]$ is the  Laplace operator,  G\"{u}nter~\cite[Appendix, \S\ IV, Thm.~3]{Gu67} has 
proved that   $W[\partial\Omega ,{\mathbf{a}},S_{{\mathbf{a}}},\cdot]$ is bounded from $C^{m-2,\alpha}(\partial\Omega)$ to $C^{m-1,\alpha'}(\partial\Omega)$ for $\alpha'\in]0,\alpha[$.

In case $n\geq 2$, $\alpha\in]0,1]$,  O.~Chkadua \cite{Chk23} has pointed out  that one could exploit Kupradze,  Gegelia,  Basheleishvili and 
 Burchuladze~\cite[Chap.~IV, Sect.~2, Thm 2.9, Chap. IV, Sect.~3, Theorems 3.26 and  3.28]{KuGeBaBu79} and  prove that if $\Omega$ is of class $C^{m,\alpha}$, then $W[\partial\Omega ,{\mathbf{a}},S_{{\mathbf{a}}},\cdot]$ is bounded from $C^{m-1,\alpha'}(\partial\Omega)$ to $C^{m,\alpha'}(\partial\Omega)$ for $\alpha'\in]0,\alpha[$.\par

In case  $n=3$, $m\geq 2$, $\alpha\in]0,1[$ and $\Omega$ is of class $C^{m,\alpha}$  and if  $P[{\mathbf{a}},D]$ is the  Helmholtz operator, Kirsch~\cite[Thm.~3.3 (a)]{Ki89} has 
developed previous work of  G\"{u}nter~\cite{Gu67}, Mikhlin~\cite{Mik70} and Colton and Kress~\cite{CoKr83} 
and has proved that  the operator  $W_\Omega[{\mathbf{a}},S_{{\mathbf{a}}}   ,\cdot]$ is  bounded from $C^{m-1,\alpha}(\partial\Omega)$ to $C^{m,\alpha}(\partial\Omega)$. 
 
 von Wahl~\cite{vo90} has considered the case of Sobolev spaces and has proved that if $\Omega$ is of class $C^{\infty}$ and if 
 $S_{{\mathbf{a}}}$ is the fundamental solution of the Laplace operator,  then  the double layer improves the regularity of one unit on the boundary. Then Heinemann~\cite{He92} has developed the ideas of  von Wahl in the frame of Schauder spaces and has proved that
 if $\Omega$ is of class $C^{m+5}$ and if 
 $S_{{\mathbf{a}}}$ is the fundamental solution of the Laplace operator,  then  the double layer improves the regularity of one unit on the boundary, \textit{i.e.}, 
 $W_\Omega[{\mathbf{a}},S_{{\mathbf{a}}}   ,\cdot]$ is linear and continuous from $C^{m,\alpha}(\partial\Omega)$ to $C^{m+1,\alpha}(\partial\Omega)$.

Maz'ya and  Shaposhnikova~\cite{MaSh05}
 have proved that $W_\Omega[{\mathbf{a}},S_{{\mathbf{a}}}   ,\cdot]$ is continuous in fractional Sobolev spaces under sharp   regularity assumptions on the boundary and if  $P[{\mathbf{a}},D]$ is the  Laplace operator. 
 
Dondi and the author \cite{DoLa17}  have proved that  if $m\geq 2$ and $\Omega$ is of class $C^{m,\alpha}$ with $\alpha\in]0,1[$, then the double layer potential $W_\Omega[{\mathbf{a}},S_{{\mathbf{a}}}   ,\cdot]$ associated to the fundamental solution of a {\em nonhomogeneous} second order elliptic differential operator with constant coefficients
 is  bounded from $C^{m,\beta}(\partial\Omega)$ to $C^{m,\alpha}(\partial\Omega)$
 for all $\beta\in]0,\alpha]$. For corresponding results for the fundamental solution of the heat equation, we refer to the  author and Luzzini   \cite{LaLu17}, \cite{LaLu18}  and references therein.

  By exploiting   a formula for the tangential derivatives of the double layer potential that involves some auxiliary integral operators  of \cite[Thm.~9.1]{DoLa17}, which generalizes the corresponding formula of 
 Hofmann, Mitrea and Taylor~\cite[(6.2.6)]{HoMitTa10} for homogeneous operators and once more by exploiting the abstract result \cite[Prop.~6.3]{La22a} in metric spaces, that generalizes previous work of Gatto \cite{Gat09}, one can prove the following.   
 For the (classical) definition of the generalized Schauder space    $C^{m,\omega_1(\cdot)}(\partial\Omega)$ we refer for example to \cite[\S 2]{DoLa17}.
\begin{theorem}
\label{wreg} 
Let ${\mathbf{a}}$ be as in (\ref{introd0}), (\ref{ellip}), (\ref{symr}). Let $S_{ {\mathbf{a}} }$ be a fundamental solution of $P[{\mathbf{a}},D]$. 
 Let $\alpha\in]0,1]$. Let $m\in{\mathbb{N}}$, $m\geq 2$. 
 Let $\Omega$ be a bounded open subset of ${\mathbb{R}}^{n}$ of class $C^{m,\alpha}$. Then the following statements hold.
 \begin{enumerate}
\item[(i)] If $\alpha\in]0,1[$, then $W_\Omega[{\mathbf{a}},S_{{\mathbf{a}}}   ,\cdot]$
 is linear and continuous from $C^{m-1,\alpha}(\partial\Omega)$ to $C^{m,\alpha}(\partial\Omega)$.
\item[(ii)] If $\alpha=1$, then $W_\Omega[{\mathbf{a}},S_{{\mathbf{a}}}   ,\cdot]$
 is linear and continuous from $C^{m-1,1}(\partial\Omega)$ to $C^{m,\omega_1(\cdot)}(\partial\Omega)$.
 \end{enumerate}
\end{theorem}
For a proof, we refer to \cite{La23c}.

Hence, Theorem \ref{wreg} sharpens the work of the above mentioned authors in the sense that if 
  $\Omega$ is of class $C^{m,\alpha}$ with $m\geq 2$, then the class of regularity of  the target space of $W_\Omega[{\mathbf{a}},S_{{\mathbf{a}}}   ,\cdot]$ is precisely $C^{m,\alpha}$  if  $\alpha<1$
  and is the generalized Schauder space $C^{m,\omega_1(\cdot)}$ if $\alpha=1$.
  
Moreover, Theorem \ref{wreg}  extends  
  the above mentioned result of  Kirsch~\cite{Ki89}  in the sense that Kirsch~\cite{Ki89}
   has  considered the Helmholtz operator in case $n=3$, $\alpha<1$ and Theorem \ref{wreg}  considers a general fundamental solution $S_{{\mathbf{a}}}$ with ${\mathbf{a}}$ as in (\ref{introd0}), (\ref{ellip}), (\ref{symr}), $\alpha\leq 1$   and $n\geq 2$.

  \section{An integral operator associated to the conormal derivative of a single layer potential}
  Another relevant layer potential operator associated to the analysis of boundary value problems for the differential operator $P[{\mathbf{a}},D]$ is defined by
\[
W_{\ast,\Omega}[{\mathbf{a}}, S_{ {\mathbf{a}} },\mu](x)\equiv
\int_{\partial\Omega}\mu(y)DS_{ {\mathbf{a}} }(x-y)a^{(2)}\nu(x)\,d\sigma_{y}\qquad\forall x\in\partial\Omega 
\] 
for all $\mu\in C^{0}(\partial\Omega)$, that we consider only in case $\Omega$ is at least of class $C^{1,\alpha}$ for some $\alpha\in ]0,1]$.  Now the continuity properties of 
$W_{\ast,\Omega}[{\mathbf{a}}, S_{ {\mathbf{a}} },\cdot]$ can be deduced by those of $W_{\Omega}[{\mathbf{a}}, S_{ {\mathbf{a}} },\cdot]$ via a simple formula. To do so, we set 
\begin{equation}\label{qrs0}
Q_j[g,\mu](x)
=\int_{\partial\Omega}(g(x)-g(y))\frac{\partial S_{ {\mathbf{a}} }}{\partial x_{j}}(x-y)\mu(y)\,d\sigma_{y}\quad\forall x\in \partial\Omega\,,
\end{equation}
for all $(g,\mu)\in  C^{0,1}(\partial\Omega)\times L^{\infty}(\partial\Omega)$, for all $l\in\{1,\dots,n\}$ and
\begin{equation}\label{eq:silapo}
V_\Omega[S_{ {\mathbf{a}} },\mu](x)\equiv 
\int_{\partial\Omega}S_{ {\mathbf{a}} }(x-y)\mu(y)\,d\sigma_{y}
\qquad\forall x\in\partial\Omega 
\end{equation}
for all $\mu\in C^{0,\alpha}(\partial\Omega)$. Then a simple computation shows that
\begin{eqnarray}
\label{v*regg1}
\lefteqn{
W_{\ast,\Omega}[{\mathbf{a}}, S_{ {\mathbf{a}} },\mu](x)
=\sum_{b,r=1}^{n}a_{br}
Q_b[\nu_{r},\mu](x)
}
\\ \nonumber
&&\qquad\qquad\qquad\qquad
-
W_\Omega[{\mathbf{a}}, S_{{\mathbf{a}}} ,\mu ](x)
-
V_\Omega[S_{{\mathbf{a}}}       ,(a^{(1)}\nu) \mu](x)\,,
\end{eqnarray}
for all $x\in \partial\Omega$ and for all $\mu\in C^{0}(\partial\Omega)$ (cf.~\cite[(10.1)]{DoLa17}) and we can prove the following statement.

  \begin{theorem}\label{v*regg}
 Let $n\in {\mathbb{N}}\setminus\{0,1\}$. 
 Let ${\mathbf{a}}$ be as in (\ref{introd0}), (\ref{ellip}), (\ref{symr}).  Let $S_{ {\mathbf{a}} }$ be a fundamental solution of $P[{\mathbf{a}},D]$.   Let $\alpha\in]0,1[$, $\beta\in]0,\alpha]$.\par 
  
 Let $\Omega$ be a bounded open subset of ${\mathbb{R}}^{n}$ of class $C^{1,\alpha}$. Then  the operator $W_{\ast,\Omega}[{\mathbf{a}}, S_{ {\mathbf{a}} },\cdot]$ is linear and continuous from $C^{0,\beta}(\partial\Omega)$ to $C^{0,\alpha}(\partial\Omega)$.
\end{theorem}
{\bf Proof.} Since the components of $\nu$ are of class $C^{0,\alpha}$, Dondi and the author  \cite[Thm.~8.2 (ii)]{DoLa17} implies that 
$Q_b[\nu_{r},\cdot]$ is continuous from $C^{0,\beta}(\partial\Omega)$ to $C^{0,\alpha}(\partial\Omega)$ for all $b$, $r\in\{1,\dots,n\}$.\par

By Theorems \ref{thm:dllregen}, \ref{thm:dllregenn} and by the continuity of the embedding of the target space of $W_{\Omega}[{\mathbf{a}}, S_{ {\mathbf{a}} },\cdot]$ into $C^{0,\alpha}(\partial\Omega)$ in each of the statements of Theorems \ref{thm:dllregen}, \ref{thm:dllregenn}, the operator $W_{\Omega}[{\mathbf{a}}, S_{ {\mathbf{a}} },\cdot]$ is linear and continuous from $C^{0,\beta}(\partial\Omega)$ to $C^{0,\alpha}(\partial\Omega)$.

Since the components of $\nu$ are of class $C^{0,\alpha}$, 
  the continuity of the pointwise product in H\"{o}lder spaces (cf.~\textit{e.g.}, \cite[Lem.~2.5]{DoLa17}) implies that the map from $C^{0,\beta}(\partial\Omega)$ to $C^{0,\beta}(\partial\Omega)$ that takes $\mu$ to $(a^{(1)}\nu) \mu$ is continuous.

By  \cite[Th.~7.2]{DoLa17}, 
and by the continuity of the embedding of $C^{0,\beta}(\partial\Omega)$ into $L^\infty(\partial\Omega)$, the operator $V_\Omega[ S_{{\mathbf{a}}}      
,\cdot]$ is continuous from $C^{0,\beta}(\partial\Omega)$ to $C^{0,\alpha}(\partial\Omega)$. 

Then formula (\ref{v*regg1}) implies  the validity of statement.\hfill  $\Box$ 

\vspace{\baselineskip}

Similarly, we can consider case $\alpha=1$ and prove the following. 
  \begin{theorem}\label{v*reggao}
 Let $n\in {\mathbb{N}}\setminus\{0,1\}$. 
 Let ${\mathbf{a}}$ be as in (\ref{introd0}), (\ref{ellip}), (\ref{symr}).  Let $S_{ {\mathbf{a}} }$ be a fundamental solution of $P[{\mathbf{a}},D]$.   Let  $\beta\in]0,1]$.\par 
  
 Let $\Omega$ be a bounded open subset of ${\mathbb{R}}^{n}$ of class $C^{1,1}$. Then  the operator $W_{\ast,\Omega}[{\mathbf{a}}, S_{ {\mathbf{a}} },\cdot]$ is linear and continuous from $C^{0,\beta}(\partial\Omega)$ to $C^{0,\omega_1(\cdot)}(\partial\Omega)$.
\end{theorem}
{\bf Proof.} Since the components of $\nu$ are of class $C^{0,1}$, \cite[Th.~8.2 (i)]{DoLa17} and the continuity of the embedding of 
$C^{0,\beta}(\partial\Omega)$ into $L^\infty(\partial\Omega)$
imply that 
$Q_b[\nu_{r},\cdot]$ is continuous from $C^{0,\beta}(\partial\Omega)$ to $C^{0,\omega_1(\cdot)}(\partial\Omega)$ for all $b\in\{1,\dots,n\}$.\par

By Theorem \ref{thm:dllreggen} (i), (ii) and by the continuity of the embedding of the target space of $W_{\Omega}[{\mathbf{a}}, S_{ {\mathbf{a}} },\cdot]$ into $C^{0, \omega_1(\cdot) }(\partial\Omega)$ in both  statements (i) and (ii), the operator $W_{\Omega}[{\mathbf{a}}, S_{ {\mathbf{a}} },\cdot]$ is linear and continuous from $C^{0,\beta}(\partial\Omega)$ to $C^{0, \omega_1(\cdot) }(\partial\Omega)$.

Since the components of $\nu$ are of class $C^{0,1}$, 
  the continuity of the pointwise product in H\"{o}lder spaces (cf.~\textit{e.g.}, \cite[Lem.~2.5]{DoLa17}) implies that the map from $C^{0,\beta}(\partial\Omega)$ to $C^{0,\beta}(\partial\Omega)$ that takes $\mu$ to $(a^{(1)}\nu) \mu$ is continuous.

Let $\beta'\in]0,\beta]\cap]0,1[$.  By  \cite[Th.~7.1 (i)]{DoLa17}, 
$V_\Omega[ S_{{\mathbf{a}}}      
,\cdot]$ is continuous from $C^{0,\beta'}(\partial\Omega)$ to
$C^{1,\beta'}(\partial\Omega)$, that is continuously embedded into $C^{0, \omega_1(\cdot) }(\partial\Omega)$. Then 
$V_\Omega[ S_{{\mathbf{a}}}      
,\cdot]$ is continuous from $C^{0,\beta}(\partial\Omega)$ to
$C^{0, \omega_1(\cdot) }(\partial\Omega)$. 

Then formula (\ref{v*regg1}) implies  the validity of statement.\hfill  $\Box$ 

\vspace{\baselineskip}

Finally, again by exploiting  formula (\ref{v*regg1}), 
one can prove the validity of the following statement in case $\Omega$ is at least of class $C^{2,\alpha}$ for some $\alpha\in]0,1]$. We also mention that the following statement   extends the corresponding result of Kirsch~\cite[Thm.~3.3 (b)]{Ki89} who has considered the case in which $S_{ {\mathbf{a}} }$ is the  fundamental solution of  the Helmholtz operator, $n=3$, $\alpha\in]0,1[$.
\begin{theorem}
\label{v*reggm}Let ${\mathbf{a}}$ be as in (\ref{introd0}), (\ref{ellip}), (\ref{symr}). Let $S_{ {\mathbf{a}} }$ be a fundamental solution of $P[{\mathbf{a}},D]$. 
 Let $\alpha\in]0,1]$. Let $m\in{\mathbb{N}}$, $m\geq 2$. 
Let $\Omega$ be a bounded open subset of ${\mathbb{R}}^{n}$ of class $C^{m,\alpha}$. 
Then the following statements hold.
 \begin{enumerate}
\item[(i)] If $\alpha\in]0,1[$,  then the operator $W_{\ast,\Omega}[{\mathbf{a}}, S_{ {\mathbf{a}} },\cdot]$ is linear and continuous from $C^{m-2,\alpha}(\partial\Omega)$ to $C^{m-1,\alpha}(\partial\Omega)$.

\item[(ii)] If $\alpha=1$,  then the operator $W_{\ast,\Omega}[{\mathbf{a}}, S_{ {\mathbf{a}} },\cdot]$ is linear and continuous from $C^{m-2,1}(\partial\Omega)$ to $C^{m-1,\omega_1(\cdot)}(\partial\Omega)$.
 \end{enumerate}
\end{theorem}
For a proof, we refer to \cite{La23c} .

\section{Examples of measures in which the lower Ahlfors regularity condition  and the doubling condition are violated}
\label{sec:exavio}
    
   We plan to present some perhaps known elementary examples of surfaces in ${\mathbb{R}}^3$ in which either the lower Ahlfors regularity inequality   or the doubling condition  are actually violated. We do so by considering   revolution graphs in ${\mathbb{R}}^3$ that are obtained by rotating a curve. 
To do so, we need some preliminary statement on the curve that we plan to rotate. If $U$ is a subset of ${\mathbb{R}}$ and if $f$ is a function from $U$ to ${\mathbb{R}}$, we say that $f$ is increasing  provided that
$f(\rho_1)\leq f(\rho_2)$ whenever $\rho_1$, $\rho_2\in U$ and $\rho_1<\rho_2$. Then we say that $f$ is strictly increasing provided that
$f(\rho_1)< f(\rho_2)$ whenever $\rho_1$, $\rho_2\in U$ and $\rho_1<\rho_2$.
\begin{proposition}\label{prop:cufn}
 Let $r_0\in]0,+\infty[$. Let $f$ be a continuous increasing function from $]0,r_0[$ to $]0,+\infty[$ such that
\begin{equation}\label{prop:cufn1}
 \lim_{x\to 0} f(x)=0\,.
\end{equation}
Then the following statements hold.
 \begin{enumerate}
\item[(i)] For each $r\in]0,r_0[$ there exists one and only one $x_r\in]0,r[$ such that
 \begin{equation}\label{prop:cufn2}
f(x_r)^2+x_r^2=r^2\,.
\end{equation}
\item[(ii)] If $f$ is also continuously differentiable, then the map from $]0,r_0[$ to $]0,r_0[$ that takes $r$ to $x_r$ is continuously differentiable and
\begin{equation}\label{prop:cufn3}
\frac{dx_r}{dr}=\frac{r}{f(x_r)f'(x_r)+x_r}\qquad\forall r\in]0,r_0[\,.
\end{equation}
\end{enumerate}
\end{proposition}
{\bf Proof.} (i) Since $f$ is increasing and continuous, the function $f(x)^2+x^2$ is strictly increasing and continuous in the variable $x\in]0,r_0[$. Now let $r\in ]0,r_0[$. Since
\[
\lim_{x\to 0} f(x)^2+x^2=0\,,\qquad f(r)^2+r^2>r^2\,,
\]
we conclude that there exists one and only one $x_r\in]0,r[$ such that equality (\ref{prop:cufn2}) holds true. 

(ii) Let $F(x,r)\equiv f(x)^2+x^2-r^2$ for all $(x,r)\in ]0,r_0[^2 $. By assumption, $F$ is continuously differentiable. Moreover, (i) implies that
\[
F(x_r,r)=0\qquad\forall r\in  ]0,r_0[\,.
\]
Also,
\[
\frac{\partial F}{\partial x}(x,r)=2f(x) f'(x)+2x\qquad\forall (x,r)\in ]0,r_0[^2 \,.
\]
Since $f$ and $f'$ are positive, we have $\frac{\partial F}{\partial x}(x_r,r)>0$ for all $r\in]0,r_0[$ and the Implicit Function Theorem implies the validity of (ii).\hfill  $\Box$ 

\vspace{\baselineskip}

Next we introduce a surface of revolution that is associated to a function 
$f\in C^1(]0,r_0[,]0,+\infty[)$. Namely, we set
\begin{equation}\label{eq:gammaf}
\gamma_f(\eta_1,\eta_2)\equiv 
\left\{
\begin{array}{ll}
 f(\sqrt{\eta_1^2+\eta_2^2}) &\text{if}\ (\eta_1,\eta_2)\in {\mathbb{B}}_2(0,r_0)\setminus\{(0,0)\}
\\
0&\text{if}\ (\eta_1,\eta_2)=(0,0)\,.
\end{array}
\right.
\end{equation}
and we plan to consider the area
\begin{eqnarray}\label{eq:arf}
\lefteqn{
A_f(r)\equiv m_2\left(
\left(
{\mathbb{B}}_3(0,r)\cap {\mathrm{graph}} (\gamma_f)
\right)\setminus\{(0,0,0)\}
\right)
}
\\ \nonumber
&&\qquad\qquad\qquad
=2\pi\int_0^{x_r}x\sqrt{1+(f'(x))^2}\,dx\qquad\forall r\in]0,r_0[\,.
\end{eqnarray}
We first show that under an extra condition on $f$, $\left({\mathrm{graph}} (\gamma_f)
\right)\setminus\{(0,0,0)\}$ is strongly upper $2$-Ahlfors regular with respect to $\{(0,0,0)\}$. 
\begin{proposition}\label{prop:cufnsua}
 Let $r_0\in]0,+\infty[$. Let $f\in C^1(]0,r_0[,]0,+\infty[)$ be increasing and satisfy the following limiting relations
 \begin{equation}\label{prop:cufnsua1}
\lim_{x\to 0}f(x)=0\,,\qquad \lim_{x\to 0}f'(x)=+\infty\,.
\end{equation}
Then the following statements hold
\begin{enumerate}
\item[(i)]  The function $x\sqrt{1+(f'(x))^2}$ is integrable in $x\in]0,r[$ for all $r\in]0,r_0[$. In particular, 
$A_f(r)<+\infty$  for all $r\in]0,r_0[$. 
\item[(ii)] $\left({\mathrm{graph}} (\gamma_f)
\right)\setminus\{(0,0,0)\}$ with the ordinary $2$-dimensional measure is strongly upper $2$-Ahlfors regular with respect to $\{(0,0,0)\}$.
\end{enumerate}
\end{proposition}
{\bf Proof.} $f$ satisfies all the assumptions of Proposition \ref{prop:cufn}. Moreover, 
de l'H\^{o}pital's rule implies that
 \begin{equation}\label{prop:cufnsua2}
\lim_{x\to 0}\frac{x}{f(x)}=\lim_{x\to 0}\frac{1}{f'(x)}=0\,.
\end{equation}
By the limiting relations (\ref{prop:cufnsua2}) and $\lim_{x\to 0}f'(x)=+\infty$,
there exists $r_0'\in]0,r_0[$ 
\[
\frac{x}{f(x)}<1\,,\qquad f'(x)>1\qquad\forall x\in]0,r_0'[\,.
\]
We now prove statement (i). Let $r\in]0,r_0[$. Since
 \[
 \lim_{x\to0}\frac{x\sqrt{1+(f'(x))^2}}{xf'(x)}=\lim_{x\to0}\sqrt{(f'(x))^{-2}+1}=1\,,
 \]
 the function $x\sqrt{1+(f'(x))^2}$ is integrable in $x\in]0,r[$ if and only if $xf'(x)$ is integrable in $]0,r[$. Since
 \[
 \int_0^{r}xf'(x)\,dx=\lim_{\epsilon\to 0}\left[xf(x)\right]_{x=\epsilon}^{x=r}-\int_0^{r}f(x)\,dx\leq f(r)r<+\infty, 
 \]
 the function $x\sqrt{1+(f'(x))^2}$ is integrable in $x\in]0,r[$ and thus statement (i) holds true.
We now prove statement (ii). We note that
\[
m_2\left(
({\mathbb{B}}_3(x,r_2)\setminus {\mathbb{B}}_3(x,r_1))\cap Y
\right)
=\int_{r_1}^{r_2}\frac{d}{dr}A_f(r)\,dr
\]
for all $r_1,r_2\in ]0,r_0'[$. Thus it suffices to show that
\begin{equation}\label{prop:cufnsua3}
\sup_{0<r<r_0'}r^{-1}\frac{d}{dr}A_f(r)<+\infty\,.
\end{equation}
Since $x_r\in]0,r[$ for all $r\in]0,r_0'[$, Proposition \ref{prop:cufn} (ii) implies that
\begin{eqnarray*}
\lefteqn{
\frac{d}{dr}A_f(r)=
2\pi x_r\sqrt{1+(f'(x_r))^2}\frac{dx_r}{dr}
}
\\ \nonumber
&&\qquad
=2\pi x_r\sqrt{1+(f'(x_r))^2}\frac{r}{f(x_r)f'(x_r)+x_r}
\\ \nonumber
&&\qquad
=2\pi \frac{x_r}{f(x_r)}\sqrt{(f'(x_r))^{-2}+1}\frac{r}{  1+(x_r/f(x_r))(f'(x_r))^{-1}}
\\ \nonumber
&&\qquad
\leq 2\pi \sqrt{1+1}r\qquad\forall r\in]0,r_0'[
\end{eqnarray*}
and thus inequality (\ref{prop:cufnsua3}) holds true and the proof of (ii) is complete.
\hfill  $\Box$ 

\vspace{\baselineskip}

Next we prove that by formulating some extra assumption on the function $f$, we can prove an apriori estimate on $x_{2r}$ and $x_r$. 
\begin{lemma}\label{lem:cufnd}
 Let $r_0\in]0,+\infty[$. Let $f\in C^1(]0,r_0[,]0,+\infty[)$ be increasing and satisfy the following limiting relations
 \begin{equation}\label{lem:cufnd1}
\lim_{x\to 0}f(x)=0\,,\qquad \lim_{x\to 0}\frac{x}{f(x)}=0\,.
\end{equation}
Let $x_r$ be as in Proposition \ref{prop:cufn}  for each $r\in]0,r_0[$. Then
 \begin{equation}\label{lem:cufnd2}
 \lim_{r\to 0}\frac{f(x_{2r})}{f(x_r)}=2
 \end{equation}
\end{lemma}
{\bf Proof.} By the definition of $x_{2r}$, $x_{r}$, we have
\[
\frac{f(x_{2r})^2+x_{2r}^2}{f(x_{r})^2+x_{r}^2}=\frac{(2r)^2}{r^2} 
\]
and accordingly
\[
\left(\frac{f(x_{2r})}{f(x_r)}\right)^2=4\frac{1+\left(x_r/f(x_r)\right)^2}{1+\left(x_{2r}/f(x_{2r})\right)^2}
\qquad\forall r\in]0,r_0/2[\,.
\]
Since $ \lim_{r\to 0}\left(x_r/f(x_r)\right)=0= \lim_{r\to 0}\left(x_{2r}/f(x_{2r})\right)$, then  the limiting relation (\ref{lem:cufnd2}) holds true.\hfill  $\Box$ 

\vspace{\baselineskip}

Next  we plan to prove a formula in order to compute 
\[
\lim_{r\to 0}\frac{A_f(2r)}{A_f(r)} \,.
\]
\begin{proposition}\label{prop:cufnnr}
 Let $r_0\in]0,+\infty[$. Let $f\in C^1(]0,r_0[,]0,+\infty[)$ be increasing and satisfy the following limiting relations
 \begin{equation}\label{prop:cufnnr1}
\lim_{x\to 0}f(x)=0\,,\qquad \lim_{x\to 0}f'(x)=+\infty\,.
\end{equation}
Let $x_r$ be as in Proposition \ref{prop:cufn} for each $r\in]0,r_0[$. If 
\[
l\equiv \lim_{r\to 0}x_{2r}/x_r\quad\text{exists\ in}\ [0,+\infty]\,,
\]
then
 \begin{equation}\label{prop:cufnnr2}
\lim_{r\to 0}\frac{A_f(2r)}{A_f(r)}=2l\,,
\end{equation}
where we understand that $2l=+\infty$ if $l=+\infty$. 
\end{proposition}
{\bf Proof.} By de l'H\^{o}pital's rule, we have
 \begin{equation}\label{prop:cufnnr3}
\lim_{x\to 0}\frac{x}{f(x)}=\lim_{x\to 0}\frac{1}{f'(x)}=0 
\end{equation}
and thus the assumptions of both Proposition \ref{prop:cufn} and Lemma \ref{lem:cufnd} are satisfied. In particular, 
\begin{equation}\label{prop:cufnnr4}
\lim_{r\to 0}\frac{f(x_{2r})}{f(x_{r})}=2\,.
\end{equation}
Moreover,  there exists $r_0'\in ]0,r_0/2[$ such that   $f'(x)\neq 0$ for all $x\in 
 ]0,2r_0'[$. By Proposition \ref{prop:cufnsua} (i),  the function $x\sqrt{1+(f'(x))^2}$ is integrable in $x\in]0,r_0/2[$. Then
 both the numerator and the denominator  
 of the fraction
\[
\frac{A_f(2r)}{A_f(r)}=\frac{\int_0^{x_{2r}}
x\sqrt{1+(f'(x))^2}\,dx}{ \int_0^{x_{r}} x\sqrt{1+(f'(x))^2}\,dx }
\]
tend to $0$ as $r$ tends to $0$. By de l'H\^{o}pital's rule, the limit of $\frac{A_f(2r)}{A_f(r)}$ as $r$ tends to $0$ exists provided that the limit of the following ratio
\begin{eqnarray*}
\lefteqn{
\frac{ 
x_{2r}\sqrt{1+(f'(x_{2r}))^2} \frac{dx_{2r}}{dr}}{   x_{r}\sqrt{1+(f'(x_{r}))^2}   
\frac{dx_{r}}{dr}}
}
\\ \nonumber
&&\qquad
=\frac{x_{2r}}{x_r}\frac{f'(x_{2r})}{f'(x_{r})}
\frac{
\sqrt{(f'(x_{2r}))^{-2}+1}
}
{
\sqrt{(f'(x_{r}))^{-2}+1}
}
\frac{
\frac{2r}{f(x_{2r})f'(x_{2r})+x_{2r}}2
}{
\frac{r}{f(x_{r})f'(x_{r})+x_{r}}
}
\\ \nonumber
&&\qquad
=4\frac{x_{2r}}{x_r}
\frac{
\frac{1}{f(x_{2r})}
}{
\frac{1}{f(x_{r})}
}\frac{
\sqrt{(f'(x_{2r}))^{-2}+1}
}
{
\sqrt{(f'(x_{r}))^{-2}+1}
}
\frac{
\frac{1}{1+(f'(x_{2r}))^{-1}(x_{2r}/f(x_{2r}) }
}{
\frac{1}{1+(f'(x_{r}))^{-1}(x_{r}/f(x_{r})}
}\quad\forall r\in]0,r_0'[
\end{eqnarray*}
exists as $r$ tends to $0$, and if such a limit exists, the two limits are equal. Then the limiting relations (\ref{prop:cufnnr1}), 
(\ref{prop:cufnnr3}), (\ref{prop:cufnnr4}) imply that
\[
\lim_{r\to 0}\frac{A_f(2r)}{A_f(r)}=4\lim_{r\to 0}\frac{x_{2r}}{x_r}\frac{1}{2}=2l\,,
\]
where we understand that $2l=+\infty$ if $l=+\infty$.\hfill  $\Box$ 

\vspace{\baselineskip}

We are now ready to present the following example of a graph of a function such that the ordinary $2$-dimensional measure fails to satisfy the doubling condition  at a point and that accordingly cannot satisfy a $2$-Ahlfors regularity condition with respect to the set that contains that point. 

\begin{example}\label{esem:badloca}
 Let $r_0\in]0,1[$. Let 
 \[
 f(x)\equiv \frac{1}{|\log x|}\qquad\forall x\in]0,r_0[\,.
 \]
 Let $\gamma_f$ be as in (\ref{eq:gammaf}). Let
 \[
Y\equiv{\mathrm{graph}}(\gamma_f)\setminus\{(0,0,0)\}
\]
 be endowed with the ordinary $2$-dimensional measure. Then
 \begin{enumerate}
\item[(i)] $Y$ is strongly upper $2$-Ahlfors regular with respect to $\{(0,0,0)\}$.

\item[(ii)] The ordinary $2$-dimensional measure $m_2$ on 
\[
{\mathrm{graph}}(\gamma_f)\setminus\{(0,0,0)\}
\]
does not satisfy the doubling condition with respect to the set $\{(0,0,0)\}$. More precisely,
 \[
 \lim_{r\to0}\frac{A_f(2r)}{A_f(r)}=+\infty 
 \]
 (cf.~(\ref{eq:arf})).
 \end{enumerate}
\end{example}
{\bf Proof.}   Since $f'(x)=\frac{1}{x\log^2x}$ for all $x\in]0,r_0[$, $f$ satisfies all the assumptions of Propositions \ref{prop:cufn}, \ref{prop:cufnsua} and thus statement (i) holds true. 

(ii) Since $f$ satisfies the assumptions of Propositions \ref{prop:cufn}, 
Proposition \ref{prop:cufnnr}, it suffices to show that
\begin{equation}\label{esem:badloca1}
\lim_{r\to 0}\frac{x_{2r}}{x_r}=+\infty\,.
\end{equation}
By Lemma \ref{lem:cufnd} there exists $r_1\in]0,r_0/2[$ such that
\[
\frac{3}{2}<\frac{ \frac{1}{|\log x_{2r}|}  }{ \frac{1}{|\log x_r|} }<\frac{5}{2}\qquad\forall r\in]0,r_1[\,,
\]
\textit{i.e.},
\[
-\frac{3}{2}\log(x_{2r})<-\log (x_r)<-\frac{5}{2}\log(x_{2r})\qquad\forall r\in]0,r_1[\,,
\]
or equivalently
\[
 x_{2r}^{-\frac{3}{2}}<x_r^{-1}<x_{2r}^{-\frac{5}{2}}
\qquad\forall r\in]0,r_1[\,.
\]
Then we have
\[
\frac{x_{2r}}{x_r}\geq  x_{2r} x_{2r}^{-\frac{3}{2}}
=x_{2r}^{-\frac{1}{2}}
\geq (2r)^{-\frac{1}{2}}
\qquad\forall r\in]0,r_1[ 
\]
and thus the limiting relation (\ref{esem:badloca1}) holds true.
\hfill  $\Box$ 

\vspace{\baselineskip}

\appendix
    
\section{Appendix: Conditions of upper Ahlfors regularity on subsets of ${\mathbb{R}}^n$ that are local Lipschitz graphs}
\label{appendix}

We first say what we mean by a subset of ${\mathbb{R}}^n$ that is a local  graph of a continuous function. 
 \begin{definition}
 \label{prelim.cocylind}  Let $n\in {\mathbb{N}}\setminus\{0,1\}$.
 Let $S$ be a subset of ${\mathbb{R}}^{n}$. Let $p\in S$, $R \in O_{n}({\mathbb{R}})$, $r$, $\delta\in ]0,+\infty[$. We say that the set
 \[
 C(p,R,r,\delta)\equiv p+ R^{t}({\mathbb{B}}_{n-1}(0,r)\times]-\delta,\delta[    ) 
 \]
 is a coordinate cylinder for $S$ around $p$, provided that the intersection
 \[
 R(S-p) \cap ({\mathbb{B}}_{n-1}(0,r)\times]-\delta,\delta[    )
 \]
 is the graph of a continuous function $\gamma$ from ${\mathbb{B}}_{n-1}(0,r)$ to $]-\delta ,\delta [$,  which vanishes at $0$ and such that $|\gamma(\eta)|<\delta/2$ for all $\eta\in {\mathbb{B}}_{n-1}(0,r)$, {\textit{i.e.}}, provided that there exists $\gamma\in C^{0}({\mathbb{B}}_{n-1}(0,r), ]-\delta ,\delta [)$ such that 
 \begin{eqnarray}
\label{prelim.cocylind1}
\lefteqn{R(S-p )\cap ({\mathbb{B}}_{n-1}(0,r)\times ]-\delta,\delta[)
}
\\   \nonumber
&&\qquad\ 
=\left\{
(\eta,y)\in 
{\mathbb{B}}_{n-1}(0,r)\times ]-\delta,\delta[:\, y=\gamma(\eta)
\right\}
\equiv{\mathrm{graph}} (\gamma) 
\,,
\\   \nonumber
&&|\gamma(\eta)|<\delta/2\qquad\forall \eta\in {\mathbb{B}}_{n-1}(0,r)\,,\qquad
\gamma(0)=0\,.
\end{eqnarray}
\end{definition}
Given a coordinate cylinder $C(p,R,r,\delta)$ for $S$ around $p$, the corresponding function $\gamma$ is uniquely determined. Indeed, if $\eta\in {\mathbb{B}}_{n-1}(0,r)$, then $\gamma(\eta)$ is the unique element $y$ of $]-\delta,\delta[$ such that
\[
(\eta,y)\in R(S-p )\cap ({\mathbb{B}}_{n-1}(0,r)\times ]-\delta,\delta[)\,.
\]
We say that $\gamma$ is the function that represents $S$ in the coordinate cylinder $C(p,R,r,\delta)$ as a graph and that the function $\psi_{p}$ from ${\mathbb{B}}_{n-1}(0,r)$ to ${\mathbb{R}}^{n}$ defined by 
\begin{equation}
\label{prelim.cocylind3}
\psi_{p}(\eta)\equiv p+R^{t}\left(
\begin{array}{c}
\eta
\\
\gamma(\eta)
\end{array}
\right)\qquad\forall \eta\in {\mathbb{B}}_{n-1}(0,r)\,,
\end{equation}
 is the parametrization of $S$ around $p$ in the coordinate cylinder $C(p,R,r,\delta)$.

Since the continuous function $\gamma$ induces the homeomorphism $(\cdot,\gamma(\cdot))$ from its domain ${\mathbb{B}}_{n-1}(0,r)$ onto its graph ${\mathrm{graph}} (\gamma)$,    the map $\psi_{p}$ is a homeomorphism from ${\mathbb{B}}_{n-1}(0,r)$ onto $\psi_{p}({\mathbb{B}}_{n-1}(0,r))=S
\cap C(p,R,r,\delta)$. 

It is sometimes useful to know that by shrinking $r$ we still obtain a    coordinate cylinder around the point $p$. More precisely, we have the following.
\begin{remark}\label{prelim.smcyl}
If $C(p,R,r,\delta)$ is a coordinate cylinder around the point $p$ of a subset $S$ of ${\mathbb{R}}^n$, then also $C(p,R,\rho,\delta)$ is a coordinate cylinder around the point $p$ of $S$ for each $\rho\in]0,r[$, and    the restriction $\gamma_{|{\mathbb{B}}_{n-1}(0,\rho)}$ represents $S$ in $C(p,R,\rho,\delta)$ as a graph.
\end{remark}

We also note that ${\mathrm{graph}} (\gamma) $ is easily seen to be path connected and that accordingly 
\[
S\cap C(p,R,r,\delta)
=p+R^{t}({\mathrm{graph}} (\gamma))
\]
 is path connected. Hence,  $S\cap C(p,R,r,\delta)$ is contained in at most one connected component of $S$.\par

We are now ready to introduce the following.
\begin{definition} Let $n\in {\mathbb{N}}\setminus\{0,1\}$.
 We say that a subset $S$ of ${\mathbb{R}}^{n}$ is a local graph of class $C^{0}$ provided that for every point $p\in S$, there exist $R\in O_{n}({\mathbb{R}})$ and $r,\delta\in]0,+\infty[$ such that $C(p,R,r,\delta)$ is a coordinate cylinder for $S$ around $p$. 
\end{definition}
If $S$ is a local   graph   of class $C^{0}$ and if $p\in S$, then
 possibly shrinking $r$,  we can always assume that $r<\delta$ and that the corresponding function $\gamma$,  
 which represents $S$ in the coordinate cylinder $C(p,R,r,\delta)$ as a graph, has a continuous extension to 
$\overline{{\mathbb{B}}_{n-1}(0,r)}$ (cf.   Remark \ref{prelim.smcyl}). It is also customary  to denote such extension by the same symbol $\gamma$. 

Then we say that $S$ is of class $C^{m}$ or of class $C^{m,\alpha}$
for some $m\in{\mathbb{N}}$, $\alpha\in ]0,1]$ provided that  $\gamma$ is of class $C^{m}$ or of class $C^{m,\alpha}$ for all $p\in\partial\Omega$. For the sake of brevity, we set  
\begin{equation}
\label{cm0}
C^{m,0}\equiv C^{m}\,.
\end{equation}
If $S$ is of class $C^{0,1}$, then we also say that $S$ is a local Lipschitz graph.   

Since $C^{1}(\overline{{\mathbb{B}}_{n-1}(0,r)})\leq  C^{0,\alpha}(\overline{{\mathbb{B}}_{n-1}(0,r)})$ for all $r\in]0,+\infty[$, $\alpha\in [0,1]$, a local graph of class $C^{1}$ is also of class $C^{0,\alpha}$. 

 Then we have the following statement that says that if $S$ is a compact local  graph of class $C^{0,\alpha}$, then we can make a uniform choice of the parameters $r$ and $\delta$. For a proof, one can follow line by  line  the corresponding proof of  
\cite[Lemma 10.1]{La19} for the case in which $S$ equals the boundary of a bounded open set of class $C^{0,\alpha}$. For the (classical) definition of norm in $C^{0,\alpha}(\overline{{\mathbb{B}}_{n-1}(0,r)})$ we refer for example to \cite[\S 2]{DoLa17}.
\begin{lemma}[of the uniform cylinders for local H\"{o}lder graphs]\label{lem:unif0cyl}
Let $n\in {\mathbb{N}}$, $n\geq 2$.    Let $\alpha\in[0,1]$. Let $S$ be a compact local  graph of class $C^{0,\alpha}$. Let $r_\ast$, $\delta_\ast\in]0,+\infty[$. Then there exist $r\in ]0,r_\ast[$, $\delta\in]0,\delta_\ast[$, $r<\delta$ such that for each  $x\in S$  there exists $R_x\in O_n({\mathbb{R}})$ such that $C(x,R_x,r,\delta)$ is a coordinate cylinder for $S$ around $x$ and the corresponding function $\gamma_x$ satisfies  the inequality
\[
 \sup_{x\in S}\|\gamma_x\|_{ C^{0,\alpha}(\overline{{\mathbb{B}}_{n-1}(0,r)}) }<+\infty\,.
 \]
\end{lemma}

Next we prove the following   extension of the  well-known fact that the boundary of a bounded open Lipschitz subset of ${\mathbb{R}}^n$ is upper $(n-1)$-Ahlfors regular with respect to itself. 
\begin{proposition}\label{prop:lgar}
 Let $n\in {\mathbb{N}}$, $n\geq 2$.   Let $S$ be a compact local  Lipschitz graph in ${\mathbb{R}}^n$, which we assume to be equipped with the ordinary $(n-1)$-dimensional measure $m_S$.   Then   $S$ is upper $(n-1)$-Ahlfors regular with respect to ${\mathbb{R}}^n$.
\end{proposition}
{\bf Proof.} Let $r$, $\delta\in]0,1[$ be as in Lemma \ref{lem:unif0cyl} of the uniform cylinders. Then we know that for each $\xi\in S$ there exist
$R_\xi\in O_n({\mathbb{R}})$ such that $C(\xi,R_\xi,r,\delta)$ is a coordinate cylinder for $S$ around $\xi$ and that the corresponding function $\gamma_\xi$ satisfies  the inequality
\[
 a\equiv\sup_{\xi\in S}\|\gamma_\xi\|_{ C^{0,\alpha}(\overline{{\mathbb{B}}_{n-1}(0,r)}) }<+\infty\,.
 \]
If $x\in {\mathbb{R}}^n$ and the distance ${\mathrm{dist}}\,(x,S)$ of $x$ from $S$ is less than $ r/4$, then there exists $x_S\in S$ such that $|x-x_S|= {\mathrm{dist}}\,(x,S)$. By the triangular inequality, we have
\[
{\mathbb{B}}_n(x, r/4)\subseteq  {\mathbb{B}}_n(x_S, r)\subseteq  C(x_S,R_{x_S},r,\delta) \,.
\]
In particular, there exists a unique $(\eta_x, y_x)\in {\mathbb{B}}_{n-1}(0,r)\times]-\delta,\delta[$ such that
\[
x=x_S+R_{x_S}^t(\eta_x, y_x)^t\,,
\qquad
|\eta_x|^2+|y_x|^2<(r/4)^2
\,.
\]
Thus if $\rho\in ]0, r/4[$, we have
\begin{eqnarray*}
\lefteqn{
S\cap {\mathbb{B}}_n(x,\rho)
 =C(x_S,R_{x_S},r,\delta)\cap S\cap {\mathbb{B}}_n(x,\rho)
 }
\\ \nonumber
&&\qquad
=x_S+R_{x_S}^t\bigl\{\bigr.
 (\eta,\gamma_{x_S}(\eta)):\,\eta\in {\mathbb{B}}_{n-1}(0,r)\,, 
 \\ \nonumber
&&\qquad\qquad\qquad\qquad\quad
  |\eta-\eta_x|^2+|\gamma_{x_S}(\eta)-y_x|^2<\rho^2
 \bigl.\bigr\}
 \\ \nonumber
&&\qquad
\subseteq
x_S+R_{x_S}^t\left\{
 (\eta,\gamma_{x_S}(\eta)):\,\eta\in {\mathbb{B}}_{n-1}(\eta_x,\rho)
 \right\}
\end{eqnarray*}
and thus
\begin{eqnarray*}
\lefteqn{
m_S(S\cap {\mathbb{B}}_n(x,\rho))
}
\\ \nonumber
&&\qquad
\leq \int_{{\mathbb{B}}_{n-1}(\eta_x,\rho)}\sqrt{1+|D\gamma_{x_S}(\eta)|^2}\,d\eta
\leq\omega_{n-1}\sqrt{1+a^2}\rho^{n-1}\,.
\end{eqnarray*}
On the other hand if  ${\mathrm{dist}}\,(x,S)\geq  r/4$, then we have $S\cap {\mathbb{B}}_n(x,\rho)=\emptyset$ for all  $\rho\in ]0, r/4[$ and thus
\[
m_S(S\cap {\mathbb{B}}_n(x,\rho))=0\leq \omega_{n-1}\sqrt{1+a^2}\rho^{n-1}\qquad\forall \rho\in ]0, r/4[\,.
\]
Hence, we conclude that statement (i) holds true and that we can choose $r_{{\mathbb{R}}^n,S,n-1}= r/4$, $c_{{\mathbb{R}}^n,S,n-1}=\omega_{n-1}\sqrt{1+a^2}$.\hfill  $\Box$ 

\vspace{\baselineskip}

For the strong upper Ahlfors regularity instead, we must formulate some extra assumption and we prove the following statement.
\begin{proposition}\label{prop:lgsar}
Let $n\in {\mathbb{N}}$, $n\geq 2$.  Let $S$ be a compact local  Lipschitz graph in ${\mathbb{R}}^n$, which we assume to be equipped with the ordinary $(n-1)$-dimensional measure $m_S$.   

 Assume that there exist $r,\delta\in ]0,+\infty[$ such that for each  $x\in S$  there exists $R_x\in O_n({\mathbb{R}})$ such that $C(x,R_x,r,\delta)$ is a coordinate cylinder for $S$ around $x$ and the corresponding function $\gamma_x$ satisfies  the inequalities
\begin{eqnarray}\label{prop:lgar1}
 &&a\equiv\sup_{x\in S}\|\gamma_x\|_{ C^{0,1}(\overline{{\mathbb{B}}_{n-1}(0,r)}) }<+\infty\,,
 \\ \nonumber
&&b\equiv \inf_{x\in S}{\mathrm{ess\,inf}}_{\eta\in{\mathbb{B}}_{n-1}(0,r)\setminus\{0\}}
\frac{(\eta+\gamma_x(\eta)D\gamma_x(\eta))\cdot\eta}{|\eta|^2}>0\,.
\end{eqnarray}
 Then  $S$ is strongly upper $(n-1)$-Ahlfors regular  with respect to $S$.
 \end{proposition}
{\bf Proof.}  We plan to prove the strong  upper $(n-1)$-Ahlfors regularity by estimating the first order derivative of $m_S(S\cap {\mathbb{B}}_n(x,\rho))$ with respect to $\rho$. To do so, we  fix $x\in S$ and we note that 
\begin{eqnarray*}
\lefteqn{
S\cap {\mathbb{B}}_n(x,\rho)
 =C(x,R_x,r,\delta)\cap S\cap {\mathbb{B}}_n(x,\rho)
 }
\\ \nonumber
&& 
=x+R_x^t\left\{
 (\eta,\gamma_{x}(\eta)):\,\eta\in {\mathbb{B}}_{n-1}(0,r)\,, 
 |\eta |^2+|\gamma_{x}(\eta)|^2<\rho^2
 \right\}
 \quad\forall\rho\in]0,r[\,,
 \end{eqnarray*}
and that accordingly
\[
m_S( S\cap {\mathbb{B}}_n(x,\rho) )= 
\int_{\{
\eta\in  {\mathbb{B}}_{n-1}(0,r):\,|\eta|^2+|\gamma_x(\eta)|^2<\rho^2
\}}\sqrt{1+|D\gamma_x(\eta)|^2}\,d\eta
 \,,
\]
for all $\rho\in]0,r[$. In order to estimate  $\frac{d}{d\rho}m_S( S\cap {\mathbb{B}}_n(x,\rho) )$, we plan to compute the integral in the right hand side by exploiting the theorem of integration for functions that are defined on domains that are normal with respect to the unit sphere. To do so however, we need to show that 
\[
A_\rho\equiv\{
\eta\in  {\mathbb{B}}_{n-1}(0,r):\,|\eta|^2+|\gamma_x(\eta)|^2<\rho^2
\}
\]
is star shaped with respect to $0$ for almost all directions of $\partial {\mathbb{B}}_{n-1}(0,1)$, \textit{i.e.}, that there exists a subset $N$ of measure $0$ of $\partial {\mathbb{B}}_{n-1}(0,1)$ such that $s\eta\in A_\rho$ for all $s\in]0,1]$ and $\eta\in A_\rho$ such that $\frac{\eta}{|\eta|}\in \partial {\mathbb{B}}_{n-1}(0,1)\setminus N$. Since $\gamma_x$ is Lipschitz continuous, the Rademacher Theorem implies that there exists
a subset $E$ of measure zero of ${\mathbb{B}}_{n-1}(0,r)$ such that $\gamma_x$ is differentiable at all points of ${\mathbb{B}}_{n-1}(0,r)\setminus E$. Then  
by applying the Theorem of integration on the spheres to the characteristic function of $E$, we can infer the existence of a subset $N$ of measure $0$ of $\partial {\mathbb{B}}_{n-1}(0,1)$ such that if $u\in \partial {\mathbb{B}}_{n-1}(0,1)\setminus N$, then $\gamma_x$ is differentiable at $su$ for almost all $s\in]0,r[$. 

Next we plan to show that for each $u\in \partial {\mathbb{B}}_{n-1}(0,1)\setminus N$ and $\rho\in]0,r[$, there exists one and only one  $r(\rho,u)\in  ]0,r[$ such that 
\[
(r(\rho,u)u,\gamma_x(r(\rho,u)u))\in \partial {\mathbb{B}}_{n}(0,\rho)\,,
\]
 \textit{i.e.}, such that
\[
|r(\rho,u)u|^2+|\gamma_x(r(\rho,u)u)|^2=\rho^2\,.
\]
To do so,  we set  
\[
G(\rho,\eta,s)\equiv |s\eta|^2+|\gamma_x(s\eta)|^2-\rho^2\quad\forall (\rho,\eta,s)\in ]0,r[\times (\partial {\mathbb{B}}_{n-1}(0,1)\setminus N)\times]0,r[
\,.
\]
 If we fix $(\rho,\eta)\in]0,r[\times (\partial {\mathbb{B}}_{n-1}(0,1)\setminus N)$, the function $G(\rho,\eta,\cdot)$ is differentiable for almost all $s\in]0,r[$ and we have
\begin{eqnarray}\label{prop:lgar2}
\lefteqn{
\frac{\partial G}{\partial s}(\rho,\eta,s)=2s\eta\cdot\eta+2\gamma_x(s\eta)D\gamma_x(s\eta)\eta
}
\\ \nonumber
&&\qquad
=\frac{2}{s}\left[
(s\eta)+\gamma_x(s\eta)D\gamma_x(s\eta)
\right]\cdot (s\eta)\geq \frac{2}{s}b|s\eta|^2=2bs|\eta|^2>0
\end{eqnarray}
for  almost  all $\ s\in]0,r[$. 
Since $G(\rho,\eta,\cdot)$ is Lipschitz continuous in $[0,r]$, 
\[
G(\rho,\eta,0)=-\rho^2<0\,,\qquad  G(\rho,\eta,\rho)\geq 0\,,
\]
 we conclude that $G(\rho,\eta,\cdot)$ is strictly increasing in $[0,r[$ and that there exists one and only one $s\in ]0,\rho]$ such that $G(\rho,\eta,s)=0$ and we set $s\equiv  r(\rho,\eta)$. 
We also note that
\begin{equation}\label{prop:lgar3}
|s\eta|^2+|\gamma_x(s\eta)|^2<\rho^2\  \forall s\in ]0,r(\rho,\eta)[\,,
\ \ 
|s\eta|^2+|\gamma_x(s\eta)|^2>\rho^2\  \forall s\in ]r(\rho,\eta),r[\,.
\end{equation}
Next we turn to show that if $\eta\in  (\partial {\mathbb{B}}_{n-1}(0,1))\setminus N$, then the function $r(\cdot,\eta)$ is Lipschitz continuous. Let $\rho_1$, $\rho_2\in ]0,r[$. There is no loss of generality in assuming that $\rho_1<\rho_2$. For the sake of brevity, we set $\varsigma_j\equiv r(\rho_j,\eta)$ for $j\in\{1,2\}$. Since $G(\rho_1,\eta,r(\rho_1,\eta))=0=G(\rho_2,\eta,r(\rho_2,\eta))$, we have
\begin{eqnarray*}
\lefteqn{
\rho_2^2-\rho_1^2= |\varsigma_2\eta|^2+\gamma_x^2(\varsigma_2\eta)
-|\varsigma_1\eta|^2-\gamma_x^2(\varsigma_1\eta)
}
\\ \nonumber
&&\qquad
=\int_0^1\frac{\partial}{\partial s}_{|s=\varsigma_1+t(\varsigma_2-\varsigma_1)}\left\{
|s\eta|^2+\gamma_x^2(s\eta)
\right\}\,dt(\varsigma_2-\varsigma_1) 
\\ \nonumber
&&\qquad
=\int_0^1\left\{
2s\eta\cdot\eta+2\gamma_x(s\eta)D\gamma_x(s\eta)\cdot\eta
\right\}_{|s=\varsigma_1+t(\varsigma_2-\varsigma_1)}\,dt(\varsigma_2-\varsigma_1)
\\ \nonumber
&&\qquad\geq 2b|\eta|^2\int_0^1\varsigma_1+t(\varsigma_2-\varsigma_1)\,dt(\varsigma_2-\varsigma_1)
\geq 2b|\eta|^2\frac{\varsigma_1+\varsigma_2}{2}(\varsigma_2-\varsigma_1) 
\end{eqnarray*}
(cf.~(\ref{prop:lgar2})). Now by the equalities
$
|\varsigma_j\eta|^2+\gamma_x^2(\varsigma_j\eta)=\rho_j^2$ for $ j\in\{1,2\}$, 
we obtain
\[
\rho_j^2\leq \varsigma_j^2+{\mathrm{Lip}}^2(\gamma_x)\varsigma_j^2=\varsigma_j^2(1+{\mathrm{Lip}}^2(\gamma_x))^2
\qquad\forall j\in\{1,2\}
\]
and thus
\[
\varsigma_2-\varsigma_1\leq\frac{1}{b}\frac{\rho_1+\rho_2}{\varsigma_1+\varsigma_2}(\rho_2-\rho_1)\leq\frac{\sqrt{1+{\mathrm{Lip}}^2(\gamma_x)}}{b}
(\rho_2-\rho_1)
\leq\frac{\sqrt{1+a^2}}{b}
(\rho_2-\rho_1)
\]
and accordingly, $r(\cdot,\eta)$ is Lipschitz continuous with Lipschitz constant less or equal to $\frac{\sqrt{1+a^2}}{b}$, which is independent of the choice of $\eta$ in $ (\partial {\mathbb{B}}_{n-1}(0,1))\setminus N$. Then by integrating on the spheres, we have
\begin{eqnarray*}
\lefteqn{
m_S( S\cap {\mathbb{B}}_{n-1}(x,\rho) )
}
\\ \nonumber
&&\qquad
= \int_{ \partial {\mathbb{B}}_{n-1}(0,1) }\int_0^{r(\rho,\eta)}
\sqrt{1+|D\gamma_x(s \eta)|^2}
s^{n-2}\, ds\,d\sigma_\eta
 \quad\forall\rho\in]0,r[\,.
\end{eqnarray*}
Since $r(\cdot,\eta)$ is Lipschitz continuous with a constant that is independent of
$\eta$  in $ (\partial {\mathbb{B}}_{n-1}(0,1))\setminus N$ and 
\[
|D\gamma_x(s\eta)|\leq a\qquad \text{a.a.}\  (\eta,s)\in ((\partial {\mathbb{B}}_{n-1}(0,1))\setminus N)\times]0,r[
\]
 and $r(\rho,\eta)\in]0,\rho]$ for all $(\rho,\eta)\in]0,r[\times((\partial {\mathbb{B}}_{n-1}(0,1))\setminus N)$, we conclude that $m_S( S\cap {\mathbb{B}}_{n-1}(x,\cdot))$ is Lipschitz continuous in $]0,r[$ and that 
\begin{eqnarray*}
\lefteqn{
\left|\frac{d}{d\rho}m_S( S\cap {\mathbb{B}}_{n-1}(x,\rho) )
\right|
}
\\ \nonumber
&&\qquad
=\left| \int_{ \partial {\mathbb{B}}_{n-1}(0,1) }
\sqrt{1+|D\gamma_x(r(\rho,\eta) \eta)|^2}\frac{\partial}{\partial \rho} r(\rho,\eta)r(\rho,\eta)^{n-2}\,  d\sigma_\eta\right|
\\ \nonumber
&&\qquad
\leq s_{n-1}\sqrt{1+a^2}\frac{\sqrt{1+a^2}}{b}\rho^{n-2} \quad
{\mathrm{a.a.}}\ 
\rho\in]0,r[ \,,
\end{eqnarray*}
where $s_{n-1}$ denotes the measure of $\partial{\mathbb{B}}_{n-1}(0,1)$ (cf.~(\ref{eq:snon})). Hence,
\begin{eqnarray*}
\lefteqn{
m_n( S\cap ({\mathbb{B}}_{n}(x,r_2)\setminus {\mathbb{B}}_n(x,r_1))  )
\leq 
\int_{r_1}^{r_2}
 s_{n-1}\frac{(1+a^2)}{b}\rho^{n-2}\,d\rho
 }
\\
&&\qquad\qquad\qquad\qquad\qquad\qquad\qquad
\leq\frac{s_{n-1}}{n-1}\frac{(1+a^2)}{b}
( r_2^{n-1}-r_1^{n-1})
\end{eqnarray*}
 for all  $x\in S$ and  $r_1,r_2\in[0,r[$ with $r_1<r_2$
 and thus we can take 
 \[
 r_{S,S,n-1}=r\,,\qquad c_{S,S,n-1}=\frac{s_{n-1}}{n-1}\frac{(1+a^2)}{b}
 \]
that are independent of the choice of $x$ in $S$  and the proof is complete. \hfill  $\Box$ 

\vspace{\baselineskip}

\begin{remark}\label{rem:coc1sa}
  If $S$   is a compact local   graph of class $C^1$, then one can follow line by line the  
proof of \cite[Lemma 2.63]{DaLaMu21}, that refers to the case in which $S$ equals the boundary of a bounded open set of class $C^{1}$ and prove prove the existence of  $r$ and $\delta$ and of $R_x$, $\gamma_x$ for all $x\in S$ such that 
\[
\sup_{x\in S}\|\gamma_x\|_{ C^{1}(\overline{{\mathbb{B}}_{n-1}(0,r)}) }<+\infty\,,
\qquad
\sup_{x\in S}\sup_{{\mathbb{B}}_{n-1}(0,r)}|D\gamma_x|<1/2\, 
\]
and deduce the validity of the conditions in (\ref{prop:lgar1}) of Proposition  \ref{prop:lgar} by means of the elementary inequality
\[
(\eta+\gamma_x(\eta)D\gamma_x(\eta))\cdot\eta\geq
\eta\cdot\eta-|\eta|^2\left(\sup_{{\mathbb{B}}_{n-1}(0,r)}|D\gamma_x|\right)^2\geq\frac{3}{4}|\eta|^2
\]
for all $\eta\in {\mathbb{B}}_{n-1}(0,r)$. Hence,   a compact local   graph $S$ of class $C^1$   is strongly upper $(n-1)$-Ahlfors regular (with respect to $S$). 
 \end{remark}

 \noindent
{\bf Acknowledgement} The author  acknowledges  the support of the Research  
Project GNAMPA-INdAM   $\text{CUP}\_$E53C22001930001 `Operatori differenziali e integrali in geometria spettrale' and is indebted to Prof.~Otari Chkadua and Prof.~David Natroshvili for a number of references, to Prof.~Joan Verdera for references \cite{Ve13} and \cite{Ve23} and to Prof.~Sergiy Plaksa for the references on the moduli of continuity  of the Cauchy integral of section \ref{dolalipc1}.

 \par\egroup

\end{document}